\documentclass{amsart}

\calclayout

\renewcommand{\baselinestretch}{1.2} 

\usepackage{graphicx} 
\usepackage{enumitem}

\usepackage[dvipsnames]{xcolor}
\usepackage{float}
\usepackage{amssymb}
\usepackage{stmaryrd}
\usepackage{amsthm}
\usepackage[utf8]{inputenc}
\usepackage[sort&compress, comma, square, numbers]{natbib}
\usepackage{upgreek}
\usepackage[bottom]{footmisc}
\usepackage[pdfpagelabels=false]{hyperref}
\hypersetup{
	colorlinks,
	linkcolor={blue!70!black},
	citecolor={blue!70!black},
	urlcolor={blue!70!black}
}
\usepackage{ulem}
\usepackage{cleveref}
\usepackage{mathtools}
\AddToHook{cmd/appendix/before}{%
    \crefalias{section}{appendix}%
    \crefalias{subsection}{appendix}%
}

\usepackage{tabularx}

\makeatletter
\renewcommand\@makefnmark{\hbox{\@textsuperscript{\normalfont\@thefnmark}}}
\makeatother

\usepackage{listings}
\usepackage{xcolor}

\lstdefinestyle{notebook}{
    backgroundcolor=\color{gray!10},
    basicstyle=\ttfamily\small,
    frame=none,
    xleftmargin=6pt,
    xrightmargin=6pt,
    aboveskip=6pt,
    belowskip=6pt,
    breaklines=true,
    columns=fullflexible,
    showstringspaces=false
}

\theoremstyle{definition}

\newtheorem{conjecture}{Conjecture}

\numberwithin{theorem}{section}

\newcommand{\co}{\varphi}

\newcommand{\acc}{\text{acc}}
\newcommand{\ord}{\text{ord}_p}

\DeclareFontFamily{OT1}{pzc}{}
\DeclareFontShape{OT1}{pzc}{m}{it}{<-> s * [1.200] pzcmi7t}{}
\DeclareMathAlphabet{\mathpzc}{OT1}{pzc}{m}{it}

\newcommand{\cL}{\mathcal{L}}

\newcommand{\cN}{\mathcal{N}}

\newcommand{\cU}{\mathcal{U}}

\newcommand{\cW}{\mathcal{W}}
\newcommand{\cX}{\mathcal{X}}
\newcommand{\cY}{\mathcal{Y}}

\DeclareFontFamily{U}{bbold}{}
\DeclareFontShape{U}{bbold}{m}{n}
{  <-5.5> s*[1.05] bbold5
	<5.5-6.5> s*[1.05] bbold6
	<6.5-7.5> s*[1.05] bbold7
	<7.5-8.5> s*[1.05] bbold8
	<8.5-9.5> s*[1.05] bbold9
	<9.5-11.5> s*[1.05] bbold10
	<11.5-16> s*[1.05] bbold12
	<16-> s*[1.05] bbold17
}{}

\newcommand{\IC}{\mathbb{C}}

\newcommand{\IF}{\mathbb{F}}

\newcommand{\IL}{\mathbb{L}}
\newcommand{\IN}{\mathbb{N}}

\newcommand{\IP}{\mathbb{P}}
\newcommand{\IQ}{\mathbb{Q}}

\newcommand{\IZ}{\mathbb{Z}}

\newcommand{\mtE}{\text{E}}

\newcommand{\mtI}{\text{I}}

\newcommand{\mtU}{\text{U}}

\newcommand{\mtW}{\text{W}}

\newcommand{\defineas}{\coloneqq}
\newcommand{\asdefine}{\eqqcolon}
\newcommand{\place}[3]{\vbox to0pt{\kern-\parskip\kern-7pt
		\kern-#2truein\hbox{\kern#1truein #3}
		\vss}\nointerlineskip}
\newcommand{\capt}[3]{\parbox{#1}{\renewcommand{\baselinestretch}{1.0}
		\caption{\label{#2}\small\it #3}}}

\newcommand{\smallfrac}[2]{\frac{\scriptstyle #1}{\scriptstyle #2}}

\renewcommand{\=}{\;=\;}

\newcommand{\Tr}{\text{Tr}}

\newcommand{\SO}{\text{SO}}

\newcommand{\GL}{\text{GL}}

\newcommand{\Teich}{\text{Teich}}
\newcommand{\Frob}{\text{Frob}}

\renewcommand{\Im}{\text{Im}}

\newcommand{\diag}{\text{diag}}

\newcommand{\me}{\text{e}}

\newcommand{\ii}{\text{i}}
\newcommand{\dd}{\text{d}}

\newcommand{\Gal}{\text{Gal}}

\newcommand{\wt}[1]{\widetilde{#1}}
\newcommand{\wh}[1]{\widehat{#1}}

\usepackage{lmodern}
\usepackage{tcolorbox}
\usepackage[absolute,overlay]{textpos}

\newtcolorbox{defn}{colback=red!5!white,colframe=red!75!black}
\newtcolorbox{funcDefn}{colback=red!5!white,colframe=red!75!black}
\newtcolorbox{objDefn}{colback=blue!5!white,colframe=blue!75!black}
\newtcolorbox{optDefn}{colback=green!5!white,colframe=green!75!black}

\newtcolorbox[auto counter,number within=section]{funcDefnN}[2][]{%
	colback=red!5!white,colframe=red!75!black,fonttitle=\bfseries,title=Function ~\thetcbcounter: #2,#1}
\newtcolorbox[auto counter,number within=section]{optDefnN}[2][]{%
	colback=green!5!white,colframe=green!75!black,fonttitle=\bfseries,title=Option ~\thetcbcounter: #2,#1}
\newtcolorbox[auto counter,number within=section]{objDefnN}[2][]{%
	colback=blue!5!white,colframe=blue!75!black,fonttitle=\bfseries,title=Object ~\thetcbcounter: #2,#1}

\title[CY Differential Operators as Truncated $p$-adic Series and Zeta Functions]{Solutions of Calabi--Yau Differential Operators as Truncated $p$-adic Series and Efficient Computation of Zeta Functions}

\author[P.Kuusela]{Pyry Kuusela}
\address{
Department of Mathematical and Physical Sciences, University of Sheffield\\
S3 7RH Sheffield, United Kingdom\\
and\\
PRISMA+ Cluster of Excellence \& Mainz Institute for Theoretical Physics\\
Johannes Gutenberg-Universit\"at Mainz\\
55099 Mainz, Germany
}
\email{pyry.r.kuusela@gmail.com}

\author[M.Lathwood]{Michael Lathwood} 
\address{Gravity, Quantum Geometry, and Field Theory Unit, Okinawa Institute of Science and Technology, 1919-1 Tancha, Onna-son, Kunigami-gun, Okinawa, 904-0495, Japan}
\email{michael.lathwood@oist.jp}

\author[M.Mosso Rojas]{Miroslava Mosso Rojas}
\address{PRISMA+ Cluster of Excellence \& Mainz Institute for Theoretical Physics\\
Johannes Gutenberg-Universit\"at Mainz\\
55099 Mainz, Germany}
\email{mmossoro@uni-mainz.de}

\author[M.Stepniczka]{Michael Stepniczka}
\address{Department of Physics, Cornell University, Ithaca, New York, 14853, USA}
\email{ms3296@cornell.edu}

\date{\today}

\begin{document}

\begin{textblock*}{5cm}(13cm,2cm) 
\raggedleft
MITP-26-011
\end{textblock*}

\begin{abstract}
    Recently, a version of the deformation method developed in ref. \cite{CdlOvS21} has been used to great effect to compute the local zeta functions of Calabi--Yau threefolds by computing their periods as series with rational coefficients and using this to find a matrix representing the Frobenius action on a $p$-adic cohomology. However, this method rapidly becomes inefficient as the prime $p$ grows, due to the rational period coefficients growing quickly. In this paper, we point out that this problem can be circumvented by a simple process that we call $p$-adically truncated recurrence. This is a recurrence relation whose solutions are $p$-adic numbers modulo $p^A$ for a given $A \in \IN$ and thus grow only slowly as $p$ grows. We show that the $p$-adic accuracy $A$ can be chosen such that all $p$-adic digits which contribute to the final result are kept, and therefore we are able to obtain the correct result by using these solutions. The improvements to speed and memory usage allow for computing the local zeta functions for tens of thousands of primes on a desktop computer, and make computing local zeta functions possible even for primes of size $10^6$ to $10^7$. Previously such computations were practically possible for around 1000 first primes. We have implemented this method in a \texttt{Sage}-compatible \texttt{Python} package \texttt{PFLFunction}.

\end{abstract}
\maketitle
\tableofcontents

\section{Introduction}
\noindent Recently, Candelas, de la Ossa, and van Straten \cite{CdlOvS21} presented a version of the Dwork deformation method for computing the local Hasse--Weil zeta functions of members $X_\co$ of one-parameter families $X_\co \hookrightarrow \cX \rightarrow \IP$ of Calabi--Yau threefolds by using the associated Picard--Fuchs operator $\cL$. Their version of the deformation method, which builds on the work of Dwork and Lauder \cite{DworkZeta, LauderCounting, LauderDeformation}, is used to find the \textit{Euler factor at prime $p$} and sufficiently nice point $\varphi$ in the moduli space
\begin{align}
\label{eq:CY3 Euler Factor}
E^{(3)}_p(X_\co, T) \; \defineas \; \det(\mtI - T \, \Frob_p^* | H^3_{\text{ét}}(X_\co, \IQ_\ell))~,
\end{align}
essentially by computing a matrix $\mtU_p(\varphi)$ representing the action of Frobenius using the solutions of the Calabi--Yau type operator $\cL$.\footnote{The deformation method of ref.~\cite{CdlOEvS20} has subsequently been generalized to cover the case of multiparameter families of threefolds of Hodge type $(1,m,m,1)$ in ref. \cite{ZetaMultThreeFolds}, see also \cite{Blesse:2025dch}. For the case of Calabi--Yau fourfolds of Hodge types $(1,1,1,1,1)$ and $(1,1,2,1,1)$, see refs. \cite{JKK23, Ducker:2025wfl}. These papers also discuss the method for elliptic curves and K3 surfaces. It should be clear how the method presented in this note can be applied to these cases as well.} Given the differential operator $\cL$ of order $b$, its solutions $\varpi^i$, $i=0,\dots,b-1$ can be computed in terms of power series $f_i(\varphi) \in \IQ \llbracket \co \rrbracket$, $i=0,\dots,b-1$, usually by deriving a recurrence relation for the coefficients. What makes practical computation possible is the assumption (supported by extensive numerical evidence) that the correct result for the Euler factor can be obtained by truncating the $\co$-series to an order $M(p,\cL)$ which depends on the prime $p$ and the operator $\cL$. 

In this note, we present a simple procedure to significantly improve the computational efficiency of this method, treating it purely as an algorithmic process having as an input a single differential operator $\cL$ of Calabi--Yau type, a prime $p$ and the truncation order $M(p,\cL)$. As the output, this method produces a polynomial $R^{(b-1)}_p(\cL, T)$. Conjecturally, this polynomial is an Euler factor of a Calabi--Yau motive.\footnote{To make a stronger statement, one can further restrict to those operators which arise as Picard--Fuchs operators of known Calabi--Yau geometries, and thus have probably an associated Calabi--Yau motive.} 

The most computationally taxing aspect of the deformation method is computing power series expansions of periods to high order: in all of the cases we have studied, the order $M(p,\cL)$ seems to grow linearly with $p$. Thus, as remarked in ref.~\cite{GvS24}, the quick growth of the rational coefficients of the power series $f_i(\co)$ causes the method to be inefficient for large primes $p$ as the amount of memory needed to store these coefficients increases rapidly. In order to circumvent this issue, we present a simple process we term \textit{$p$-adically truncated recurrence}. This is a recurrence relation whose solutions are $p$-adic numbers modulo $p^A$ for a given $A \in \IN$, and which approximate the solutions of $\cL$ truncated to order $M(p,\cL)$ in the following sense: We show that given any $B \in \IN$, we can always choose $A$ such that the solutions to the $p$-adically truncated recurrence agree with the solutions of the Calabi--Yau operator modulo $p^B$. We also find the sufficient $B$ such that the Euler factor computed using the solutions to the $p$-adically truncated recurrence agrees with that computed using the exact rational coefficients. Using the $p$-adically truncated recurrence speeds up the computations and significantly decreases the memory usage compared to computation using exact rational coefficients. 

We have implemented the method of the $p$-adically truncated recurrences in a \texttt{Sage}-compatible \texttt{Python} package \texttt{PFLFunction}.\footnote{\href{https://github.com/PyryKuusela/PFLFunction/tree/main}{\texttt{https://github.com/PyryKuusela/PFLFunction}}} With the package, we are able to compute local zeta functions of Calabi--Yau threefolds for tens of thousands of primes on a desktop computer. In addition, we compute individual zeta functions for primes of size $10^6$ to $10^7$. Previously such computations were practically possible for around 1000 first primes \cite{CdlOEvS20,CdlOvS21,GvS24}. 

We conclude by giving some brief examples of such computations and highlight some applications that are now tractable with $p$-adically truncated recurrence. We compute the Euler factors at $p=2^{20}-3$ of the family of mirror quintic threefolds, cross-checking this result with that computed using the \texttt{controlledreduction} algorithm of ref.~\cite{controlledreduction}. In the mathematical physics community, the factorization of the Euler factors $E^{(3)}_p(X_\co, T)$ over $\IQ$ has been used to study the attractor mechanism of $\cN=2$ supergravity and flux vacua, among other physical systems of significant interest \cite{CdlOEvS20, Candelas:2023yrg}. We argue that with the more efficient method it is possible to study the statistics of Frobenius traces, which can probe finer properties of Galois representations relevant to physicists. As a toy example, we find singular members of a one-parameter family of K3 surfaces using Sato--Tate distributions. Finally, we revisit the computation of ref.~\cite{GvS24}, where the deformation method is utilized to identify paramodular forms from Calabi--Yau differential operators, since analogous computations provide highly non-trivial consistency checks for the Euler factors produced by our method. 

This paper is organized as follows: In \Cref{sec:review}, we review the deformation method.
In \Cref{sec:pAdicTrunc}, we describe the $p$-adically truncated recurrence, and give a simple upper bound for the sufficient $p$-adic accuracy $A$ needed to obtain the correct Euler factor. The examples and applications are presented in \Cref{sec:examples}. In \Cref{app:Integrality}, we include, for completeness, the integrality conditions satisfied by Calabi--Yau operators. \Cref{app:Rationality} gives some details on rational structures associated to Calabi--Yau type differential operators of order 4. In \Cref{app:PFLFunction}, we include documentation for the package \texttt{PFLFunction}.

\begin{table}[H]
	\renewcommand{\arraystretch}{1.35}
	\centering
	\begin{tabularx}{\textwidth}{|>{\hsize=.12\hsize\linewidth=\hsize}X|
			>{\hsize=0.78\hsize\linewidth=\hsize}X|>{\hsize=0.05\hsize\linewidth=\hsize}X|}
		\hline
		\textbf{Symbol} & \hfil \textbf{Definition/Description} & \hfil \textbf{Ref.}\\[3pt]
		\hline	\hline
		
		$\cL$
		& 
		A differential operator of Calabi--Yau type of order $b$ and degree $N$.
		&
		\eqref{eq:PF_equation}
		\\[4pt] \hline

        $\varpi^i$
		& 
		Solutions $\cL\varpi^i=0$ to the Calabi--Yau operator $\cL$.
		&
        \eqref{eq:frobenius_basis_definition}
  	    \\[4pt] \hline

        $f_i^{[M]}(\co)$
		& 
		Holomorphic power series entering the $\varpi^j$, truncated to $M$ terms.
		&
        \eqref{eq:truncated_f}
  	    \\[4pt] \hline

        $c_{i,n}$
		& 
		The coefficient of $\co^n$ in the power series expansion of $f_i(\co)$.
		&
        \eqref{eq:period_series_definition}
  	    \\[4pt] \hline

        $c^{(A)}_{i,n}$
		& 
		Solutions to the $p$-adically truncated recurrence with initial accuracy $A$.
        &
        \eqref{eq:truncated_recurrences}
  	    \\[4pt] \hline

        $\mtE(\co)$
		& 
		The (truncated) logarithm-free period matrix.
        &
        \eqref{eq:E_tilde_definition}
  	    \\[4pt] \hline

        $\mtU_p(\co)$
		& 
		The matrix (conjecturally) corresponding to the action of the inverse Frobenius.
        &
        \eqref{eq:U_matrix_definition}
  	    \\[4pt] \hline

        $\alpha_i$
		& 
		The (prime-dependent) constants determining the matrix $\mtU_p(0)$.
        &
        \eqref{eq:U(0)}
  	    \\[4pt] \hline        

    	$R^{(b-1)}_p(\cL,T)$
		& 
		The characteristic polynomial of the matrix $\mtU_p(\Teich(\co))$ associated to operator $\cL$.
		&
        \eqref{eq:Euler_factor_p}
  	    \\[4pt] \hline

        $a^{(i)}_p$
		& 
		The coefficient of $T^i$ in $R^{(b-1)}_p(\cL,T)$.
		&
        \eqref{eq:a_p_definition}
  	    \\[4pt] \hline
        
        $A$
		& 
		The initial accuracy of the $p$-adically truncated recurrence.
        &
        \S\ref{sec:pAdicTrunc}
  	    \\[4pt] \hline

        $B$
		& 
		The $p$-adic accuracy to which $\mtU_p(\co)$ is computed.
        &
        \eqref{eq:apBound}
  	    \\[4pt] \hline

        $C$
		& 
		A rational number such that the numerator of $\mtU_p(\co)$ truncates to order $\lceil Cp \rceil$ in $\co$.
        &
        \eqref{eq:trunc_order}
  	    \\[4pt] \hline
	\end{tabularx}
	\capt{\linewidth}{tab:notation}{Some quantities that are used throughout the paper with references to where they are first introduced. }
\end{table}

\section{A Brief Review of the Deformation Method}
\label{sec:review}
\noindent To provide the necessary context for the present note, we begin with a brief review of the version of Dwork's deformation method developed in ref.~\cite{CdlOvS21}. Here we focus almost exclusively on the algorithmic and computational aspects of the method, with a brief discussion of the conceptual aspects in the next section. For an in-depth discussion, we refer the interested reader to refs.~\cite{CdlOvS21,ZetaMultThreeFolds}. 

\subsection{Differential operators of Calabi--Yau type}
\label{subsec:cyoperators}
Let $\cL$ be a differential operator
\begin{align} \label{eq:PF_equation}
\cL \; \defineas \; \sum_{i=0}^{b} S_i(\co) \; \theta^i \in \IZ[\co][\theta] \qquad \text{where} \qquad  \theta \; \defineas \; \varphi \, \partial_\varphi~,
\end{align}
with $S_i \in \IZ[\co]$. We call $b$ the \textit{order} of the operator and the maximum of the degrees of the polynomials $S_i(\co)$ corresponds to the \textit{degree} $N$ of the operator. We say that $\cL$ is of \textit{Calabi--Yau type} \cite{Bogner,DucoCYOps} if the following properties are satisfied:
\begin{enumerate}[label=\Alph*]
    \item The operator $\cL$ is of \textit{Fuchsian} type, which means that $\cL$ has only \textit{regular singular points}. A singular point $\co^*$ is said to be regular if all solutions to the differential equation $\cL = 0$ grow at most as a power of the inverse distance to $\co^*$.
    \item The point $\co = 0$ is a \textit{point of maximal unipotent monodromy} (MUM-point), meaning that the indicial equation has a single solution which is an integer with multiplicity $b$. \label{cond:MUM}
    \item The operator $\cL$ is \textit{self-dual}. This means that there exists a non-zero function $ \alpha \in \mathbb{Q}(\co)$ such that $\cL \alpha = \alpha \cL^{*}$, where $\cL^{*}$ is the adjoint of the differential operator $\cL$ defined as:\footnote{If such $\alpha$ exists, it satisfies the differential equation $\partial_{\co} \, \alpha = -\frac{2}{b} \,A_{b-1}(\co) \, \alpha$.} 
    \begin{equation}
            \cL^{*} \; \defineas \; \sum_{i=0}^{b}  (-1)^{b+i} \, \partial_{\co}^{i} \, A_{i}(\co) \in \IZ[\co][\theta]~,
    \end{equation}
    where the functions $A_i(\co)$ are defined by writing $\cL$ as
    \begin{align}
        \cL \= \partial_{\co}^{b} +  \displaystyle \sum_{i=0}^{b-1} A_{i}(\co) \, \partial_{\co}^{i}~.
    \end{align}
\end{enumerate}
In addition, the operator must satisfy three integrality conditions, which we will not directly utilize, but which are included, for completeness, in \Cref{app:Integrality}. Following the convention of ref.~\cite{DucoCYOps}, two Calabi--Yau type operators are considered to be equivalent if they are related to each other by a coordinate transformation fixing the MUM-point at the origin or by multiplication with an algebraic function.
    
We concentrate on studying to Calabi--Yau operators which have a canonical rational structure, which we can identify by using the \textit{monodromy representation} of $\cL$.\footnote{Note that not all of the operators satisfying the above conditions have such a rational structure. See ref.~\cite{DucoCYOps} for an example.} Let $\IL_\IC^{\vphantom{Q}}$ denote the local system of $\cL$, $\IL_{\IC}^{*\vphantom{Q}}$ the fiber at a fixed point $\varphi_*$, and $\Sigma$ its set of singularities. Then the monodromy representation is the map
\begin{equation}
    \begin{split}
    T : \, \pi_{1} \left( \mathbb{P}^{1} \setminus \Sigma , \varphi_* \right) \,&\longrightarrow \, \mathrm{Aut}( \IL_\IC^{*\vphantom{Q}} ) \, \simeq \, \GL_{b}( \mathbb{C})~,\\
     \gamma \, &\longmapsto \, T_{\gamma}~,
    \end{split}
\end{equation}
associating each loop $\gamma$ to the corresponding monodromy $T_\gamma$ around $\gamma$. The image $\Im \left( \pi_{1} \left( \mathbb{P}^{1} \setminus \Sigma , \varphi_* \right) \right) \subseteq \GL_{b}( \mathbb{C})$ is called the \textit{monodromy group} of $\cL$. We say that $\IL_\IQ^{\vphantom{Q}}$ is a rational structure on $\IL_{\IC}^{\vphantom{Q}}$ if the fiber $\IL_\IQ^{*\vphantom{Q}}$ at point $\varphi_*$ is a $\IQ$-vector space $\IL_\IQ^{*\vphantom{Q}} \subset \IL_\IC^{*\vphantom{Q}}$ such that $\IL_\IQ^{*\vphantom{Q}} \otimes_{\IQ}^{\vphantom{Q}} \IC \raise0.15ex\hbox{$\, \, = \, \,$} \IL_\IC^{*\vphantom{Q}}$ and $\IL_\IQ^{*\vphantom{Q}}$ is preserved under the action of monodromy representation. In other words, $\IL_\IQ^{\vphantom{Q}}$ is fixed by requiring that given any basis of $\IL_\IQ^{*\vphantom{Q}}$, the monodromy matrices have rational entries.

The condition \eqref{cond:MUM} implies that $S_{i}(0) = 0$ for $i<b$ and $S_b(0) \neq 0$, so rescaling the coordinate $\co$ if necessary, we can normalize $S_b(0)=1$. Consequently, in some open disk $\co \in D \ni 0$ around the MUM-point, one can find a basis of solutions $\varpi^i \in \IQ\llbracket\varphi\rrbracket[\log \varphi]$ annihilated by $\cL$, which one can write as
    \begin{align} \label{eq:frobenius_basis_definition}
    \varpi^i \= \sum_{m=0}^i \frac{\log^m \co}{m!} f_{i-m}(\co), \quad \text{for} \quad 0 \leq i \leq b-1 ~,
    \end{align}
where $f_i(\co) \in \IQ\llbracket\co\rrbracket$ are holomorphic functions of $\co$ which we normalize so that $f_i(0) = \delta_{i,0}$. By slight abuse of language, borrowing the geometric terminology, we call these solutions \textit{periods}. We denote the Laurent coefficients of the functions $f_i(\co)$ as $c_{i,n}$ so that we can express them as 
\begin{align} \label{eq:period_series_definition}
      f_i(\co) \; \asdefine \; \sum_{n=0}^\infty c_{i,n} \, \co^n~.
\end{align}
These \textit{period coefficients} $c_{i,n}$ satisfy recurrence relations of the form
\begin{align} \label{eq:period_recurrences}
n^b c_{i,n} \= -\sum_{j = 0}^{i-1} \frac{b!}{(b-i+j)!} n^{b-i+j} c_{j,n} + \sum_{j = 0}^i \, \sum_{k = 1}^{N} \wt P_{j,n-k}(n) \; c_{j,n-k}~, 
\end{align}
where $\wt P_{i,k}(n) \in \IZ[n]$ are polynomials with integral coefficients, and $N$ is the degree of the operator $\cL$. In writing this expression, we have used that the Calabi--Yau conditions imply that $S_i(\co)$ has no constant term for $i<b$, we have normalized so that $S_b(0)=1$. The initial conditions for the recurrence follow from the normalization of the $f_i(\co)$:
\begin{align} \label{eq:period_recurrence_initial_conds}
c_{i,0}\= \delta_{i,0}~, \qquad c_{i,n} \= 0 \text{ for } n<0~. 
\end{align}

It turns out to be useful to view the functions $f_i(\co)$ in terms of a double expansion: we have the $\varphi$-adic series expansion of $f(\co)$ with Laurent coefficients $c_{i,n} \in \IQ$, which we can expand in terms of their $p$-adic digits $c_{i,n;m}$ so that
\begin{align}
\label{eq:fPadicExpansion}
    f_i(\co)  \= \sum_{n=0}^\infty \quad \sum_{m = \ord(c_{i,n})}^\infty c_{i,n;m} \quad p^m \, \co^n~.
\end{align}
The truncation of the $\varphi$-adic series is a key element of the deformation method and has been previously studied in detail for example in refs.~\cite{Candelas:2007mb,CdlOvS21,LauderCounting,LauderDeformation}. Here, we study the truncation of the $p$-adic series. 
\subsection{Computing the polynomials \texorpdfstring{$R_p^{(b-1)}(\cL,T)$}{R}}
\label{subsec:EulerFactors}
We now describe the version of the deformation method developed in ref.~\cite{CdlOvS21}. The algorithm takes as input a Calabi--Yau differential operator $\cL$, a prime $p$, and an integer $M(\cL,p)$, which denotes the order in $\varphi$ to which the series $f_i(\varphi)$ defined in eq.~\eqref{eq:period_series_definition} are truncated. As output the procedure gives a polynomial $R_p^{(b-1)}(\cL,T)$,  which is conjecturally equal to the Euler factor $E_p^{b-1}(X_\varphi,T)$ of a Calabi--Yau motive $X_\co$ whose Picard--Fuchs operator is $\cL$. We write this polynomial as
\begin{equation}
\label{eq:Euler_factor_p}
    R_p^{(b-1)}(\cL,T) \; \defineas \; \sum_{i=0}^{b} a_p^{(i)} T^i.
\end{equation}
We now describe the process of computing these coefficients $a_p^{(i)}$ starting from the truncated power series
\begin{align}
\label{eq:truncated_f}
f_i^{[M(\cL,p)]}(\varphi) \; \defineas \;  \sum_{n=0}^{M(\cL,p)} c_{i,n} \, \co^n~,
\end{align}
which enter the \textit{(truncated) logarithm-free period matrix}\footnote{In ref.~\cite{CdlOvS21} and subsequent works, this is typically denoted by $\widetilde{\mtE}(\varphi)$, while $\mtE(\co)$ denotes the full period matrix which includes the logarithmic terms. However, we will not need this latter matrix, and thus use the less cumbersome notation $\mtE(\co)$ for the logarithm-free matrix.}, that is
\begin{align} \label{eq:E_tilde_definition}
[\mtE(\varphi)]^{i,k} \; \defineas \; \left. \theta^k \varpi^i \right|_{\log \varphi \to 0} \= \sum_{j=0}^i \binom{k}{i-j} \,\theta^{k+j-i} \, f_j^{[M(\cL,p)]}(\varphi)~.
\end{align}
The coefficients $a_p^{(i)}$ will be obtained from the traces of the \textit{$\mtU$-matrix}
\begin{align} \label{eq:U_matrix_definition}
\mtU_p(\co) \; \defineas \; \mtE(\co^p)^{-1} \, \mtU_p(0) \, \mtE(\co)~.
\end{align}
Therefore, the $\mtU$-matrix evaluated at an arbitrary point is itself determined by two pieces of data: the logarithm-free period matrix $\mtE(\co)$ and the form of the $\mtU$-matrix at the MUM point $\co=0$. The functions $f_j^{[M(\cL,p)]}(\co)$ can be obtained using the recurrence of eq.~\eqref{eq:period_recurrences}, and thus also $\mtE(\co)$. 

The inversion of $\mtE(\co)$ can be carried out using the method developed in ref.~\cite{Thorne}.\footnote{The original idea of this method relies on unpublished work by Duco van Straten, while the full procedure and practical implementation were later developed in ref.~\cite{Thorne}.} Define the matrix $\mtW(\co)$ as
\begin{equation}
    \mathrm{W}(\varphi) \; \defineas \; \mtE(\co)^{T} \, \sigma \, \mtE(\co)~, \qquad [\sigma]_{kl} \defineas (-1)^{l} \delta_{k+l, b-1}~.
\end{equation}
The utility of this matrix is that it is a matrix of rational functions in $\varphi$, and thus quick to invert in practice.\footnote{See, for instance, refs.~\cite{ZetaMultThreeFolds, SorenThesis, KlemmEtAl} for explicit expressions for operators of low degree.} The matrix $\mtE(\co)^{-1}$ can be then expressed as
\begin{align} \label{eq:E_inversion}
    \mtE(\co)^{-1} \=  \left(\sigma \, \mtE(\co) \mtW(\co)^{-1} \right)^T~,
\end{align}
reducing the problem of inverting the matrix $\mtE(\co)$ is to inversion of a matrix of rational functions.

The matrix $\mtU_p(0)$ is a constant matrix of the form
\begin{align} \label{eq:U(0)}
\begin{split}
&\mtU_p(0) \defineas u_p \, \text{diag}(1,p,\dots,p^{b-1})\sum_{i=0}^{b-1} \alpha_i \epsilon^i~,\\
&\text{with} \quad \alpha_i \in \IQ_p~,~~ \alpha_0 = 1~,~~ u_p = \pm1~,~~\epsilon \defineas \begin{pmatrix}
0 & 0 & 0 & \cdots & 0 \\
1 & 0 & 0 & \cdots & 0 \\
0 & 1 & 0 & \cdots & 0 \\
\vdots & \ddots & \ddots & \ddots & \vdots \\
0 & \cdots & 0 & 1 & 0
\end{pmatrix}
\end{split}
\end{align}
where the coefficients $\alpha_i$ satisfy the relations arising from the condition
\begin{align}
\mtU_p(0) \, \sigma \, \mtU_p(0)^{T} = p^{b-1} \sigma~.
\end{align}
Conjectural expressions for the remaining independent constants $\alpha_i$ have been given in refs.~\cite{CdlOvS21,JKK23,ZetaMultThreeFolds, Ducker:2025wfl} for differentials operators with $b=3,4,5$, for which it is also conjectured that $u_p=1$ for all $p$.\footnote{To be more specific, the conjecture states that choosing $u_p=1$ for all $p$ is a consistent choice. That is, there exists a Calabi--Yau motive $X_\co/\IQ$ such that the choice $u_p=1$ produces the Euler factor corresponding to this motive. In general, there are other consistent choices for the values of $u_p$ related to the existence of quadratic twists of $X_\co$. This is natural, since the Euler factors are associated to motives defined over $\IQ$ whereas the differential operator $\cL$ is (conjecturally) associated to a family defined over $\IC$.} For example, for operators with $b=4$, it is conjectured that (subject to a judicious choice of the coordinate $\co$, see ref.~\cite{ZetaMultThreeFolds}), there is only one independent non-zero coefficient, $\alpha_3$.\footnote{This choice of coordinates is explained in details in ref.~\cite{ZetaMultThreeFolds}. However, these coordinates do not always coincide with the coordinates normalized such that $S_{b}(0)=1.$} This takes the form $\alpha_3 = K \zeta_p(3)$, where $K \in \IQ$ is a constant independent of the prime $p$ and determined by the rational structure associated to $\cL$. See \Cref{app:Rationality} for further details.

We compute the series appearing in the matrix $\mtU(\co)$ modulo $p^B$, where $B$ is the least integer such that
\begin{align}
\label{eq:apBound}
\binom{b}{\lceil b/2 \rceil} p^{\lceil b/2 \rceil (b-1)/2} \leq p^{B}~,
\end{align}
and resum the series appearing in the resulting matrix in the form
\begin{eqnarray}
\label{eq:UratMatrix}
    \mtU_p(\co) &\eqqcolon& \frac{\cU_p(\co)}{D(\co^p)} \mod p^{B} ~,     
\end{eqnarray}
where $\cU_p(\co)$ is a matrix whose entries are polynomials in $\IZ_p[\co]$ and $D(\co)$ is, in general, an algebraic function dependent on the discriminant of the operator $\cL$.\footnote{For instance, for operators with $b=4$ the denominator often takes polynomial form $D(\co) = \cY(\co)^2 \cW(\co)$, where $\cW(\co)$ vanishes exactly at F-points, that is, singular points around which the local monodromy is of finite order, and $\cY(\co)$ vanishes at K-points, where there are two pairs of equal local indices. \cite{ZetaMultThreeFolds}.\label{foot:U_denominator}}

When $D(\Teich(\varphi))^{-1}$ is defined in $\IQ_p$ and $\ord(D(\Teich(\varphi))) \leq 0$, we define the coefficients $a_{p}^{(k)}$ for $k=0,\dots,\lceil b/2 \rceil$ recursively:
\begin{equation} \label{eq:a_p_definition}
    a_p^{(0)} \defineas 1~, \qquad a_{p}^{(k)} \defineas - \frac{1}{k}\sum_{i = 1}^{k} \Tr \, [ \mtU_{p}(\Teich(\co))^i ] \, a_{p}^{(k-i)} \mod p^B ~,
\end{equation}
where the $\mtU_p(\Teich(\co))$ matrix is evaluated at the Teichm{\"u}ller lift $\Teich(\co)$ of $\co \in \mathbb{F}_p$. For instance, the first few coefficients take the form:
\begin{equation}
\label{eq:traceFormulae}
    \begin{split}
    a_{p}^{(0)} \= & 1 , \quad a_{p}^{(1)}  \= - \Tr \, [ \mtU_p(\Teich(\co))], \quad a_{p}^{(2)} \= \frac{1}{2} \left( \Tr \, [\mtU_p(\Teich(\co)) ]^2 - \Tr \, [\mtU_p(\Teich(\co))^2] \right)~.
    \end{split}
\end{equation}
The coefficients $a_p^{(k)}$ for $k=\lceil b/2 \rceil+1,\dots,b$ are computed from these using the relation
\begin{equation} \label{eq:functional_equation}
    a_{p}^{(j)} \= a_{p}^{(b-j)} \, p^{-\frac{1}{2}(b-1)(b-2j)}~.
\end{equation}

\subsection{ \texorpdfstring{$R_p^{(b-1)}(\cL,T)$}{R} and Euler factors of Calabi--Yau motives}
As we have alluded to previously, the utility of the polynomials $R_p^{(b-1)}(\cL,T)$ comes from the following conjecture: 
\begin{conjecture}
    The Calabi--Yau operator $\cL$ arises as a Picard-Fuchs operator for some Calabi--Yau motive $X_\co$. Let $p \geq 5$ and let $\co$ be such that $D(\Teich(\varphi))^{-1}$, described in eq.~\eqref{eq:UratMatrix}, is defined in $\IQ_p$ and $\ord(D(\Teich(\varphi))) \leq 0$. Then, for sufficiently large $M(\cL,p)$, the Euler factor of this motive, $E_p^{(b-1)}(X_\co,T)$, is equal to the polynomial $R^{(b-1)}_p(\cL,T)$.
\end{conjecture}
To obtain the Euler factor, the order $M(\cL,p)$ to which the series $f_i(\co)$ are computed needs to be sufficiently large so that the resummed expression \eqref{eq:UratMatrix} is exact. One can experimentally find a sufficiently large value of $M(\cL,p)$ by computing the series $f_i(\co)$ to a high order, resum the series appearing in the matrix $\mtU_p(\co)$ in the form \eqref{eq:UratMatrix}, and finding the degrees of the polynomials appearing as coefficients of $\cU_p(\co)$. Doing this, it was found in ref.~\cite{CdlOvS21} that it appears to be possible to find a prime-independent constant $C \in \IQ$ such that a sufficient $M(p,\cL)$ is given by
\begin{align}\label{eq:trunc_order}
M(p,\cL) \=  \left \lceil C p \right \rceil~.
\end{align}
Further, the value of $C$ was linked to the indices of the operator $\cL$ at infinity. Since we are mainly concerned with computing the zeta functions for larger primes, in the following \Cref{sec:pAdicTrunc} we will often assume that we consider primes $p > C$. In practice, $C$ seems to be a small rational number, at least for low degree operators $\cL$, so this is not a very restrictive requirement.\footnote{Although the simplest expressions are attained when this requirement is met, the formulae we give also generalize in an obvious manner to cases $C<p$ and the case where $C$ is dependent on the prime.}

The main reason the statement above remains conjectural is that it is not proven in general that a Calabi--Yau motive exists for every Calabi--Yau type differential operator. If one assumes this, then the only conjectural statements concern the operator-dependent constants appearing in the matrix $\mtU_p(0)$ and the validity of the resummation formula \eqref{eq:UratMatrix}. These are supported by extensive numerical evidence (see e.g. refs.~\cite{CdlOEvS20,CdlOvS21,Candelas:2023yrg,GvS24, ZetaMultThreeFolds} for various consistency checks) and have been proven in some cases \cite{Shapiro2009,Shapiro2012,Beukers2025}. For a more careful discussion of these issues, see refs.~\cite{CdlOvS21,ZetaMultThreeFolds}.

Assuming this conjecture also elucidates the above algorithm. The Euler factor can be computed as 
\begin{equation}
    E_p^{(b-1)}(X_\varphi,T) \= \det(\text{I}-T(\Frob_p^*)^{-1}|H^{b-1}(X_\varphi,\mathbb{Q}_p)),
\end{equation}
where $\Frob_p^*$ is the pullback of the Frobenius action on an appropriate $p$-adic cohomology, $H^{b-1}(X_\varphi,\mathbb{Q}_p)$. When the differential operator $\cL$ is a Picard--Fuchs operator associated to a Calabi--Yau motive, the matrix $\mtU_p(\co)$ represents this action (up to the conjectural form of $\mtU_p(0)$). In particular, the appearance of the period matrix is a consequence of the compatibility of the Frobenius action with the Gauss--Manin connection, which can be compactly stated as
\begin{align}
\nabla_{\theta} \, \Frob_p \= p \, \Frob_p \nabla_{\theta}~.
\end{align}
The form \eqref{eq:U(0)} of the matrix $\mtU(0)$ can be derived from the $\varphi \to 0$ limit of this equation, together with a compatibility of the Frobenius map with the Poincaré intersection product.

The sufficient $p$-adic accuracy $B$ in eq.~\eqref{eq:apBound} to which we compute the matrix $\mtU(\co)$ follows from the Riemann hypothesis part of the Weil conjectures, the coefficients $a_p^{(i)}$ are those of the characteristic polynomial $\det(\mtI - T \, \mtU(\co))$, and the relation \eqref{eq:functional_equation} follows from the functional equation of the zeta function.

\section{ \texorpdfstring{$p$}{p}-adically Truncated Recurrence}
\label{sec:pAdicTrunc}
\noindent To compute the truncated series $f_i^{[Cp]}(\co)$ which are needed to find the polynomials $R^{(b-1)}_p(\cL,T)$ introduced in the previous section, one can use the recurrence relations \eqref{eq:period_recurrences} which in principle determine the period coefficients $c_{i,n}$ to an arbitrary order. However, in practice, these grow rapidly, which makes storing the series memory-intensive and computations involving them very slow. In this section, we show that this problem can be circumvented to a significant degree by a procedure we term \textit{$p$-adically truncated recurrence}: We set up a recurrence whose solutions are $p$-adic integers modulo $p^A$, for a given $A \in \IN$, and which approximate the exact period coefficients $c_{i,n}$ in the sense that their difference can be made arbitrarily small in the $p$-adic norm by choosing sufficiently large $A$. We find a lower bound for $A$ such that the accuracy provided is sufficient to guarantee that the solutions can be used to reproduce the polynomial $R^{(b-1)}_p(\cL,T)$. The simple bound we give in this section is universal, meaning that it depends only on the order $b$ of the differential operator $\cL$ and the associated $C$, not on its precise form, and does not depend on the prime $p$, which is important for the scalability of the method for higher values of $p$.

\subsection{\texorpdfstring{$p$}{p}-adic truncation of solutions} \label{sect:series_truncation}
As reviewed in \Cref{subsec:EulerFactors}, the deformation method uses the truncated $\co$ series
\begin{align}
\label{eq:fWithC}
    f^{[Cp]}_i(\co)  \= \sum_{n = 0 \vphantom{\ord(c_{i;n})}}^{\left \lceil{Cp} \right \rceil} \hskip0pt \sum_{m = \ord(c_{i,n})}^\infty \hskip-10pt c_{i,n;m} \, p^m \, \co^n
\end{align}
to compute the $p$-adic expansion of $\mtU_p(\co)$ up to order $B$. Unlike the truncation of the $\varphi$ expansions, which only ever involve positive powers of $\varphi$, the $p$-adic expansion is slightly more subtle due to the fact that when $n$ grows, the coefficients $c_{i,n}$ tend to involve larger negative powers of $p$. We must therefore find the $p$-adic valuations of the coefficients of the series expansions that enter the formula \eqref{eq:U_matrix_definition} for the matrix $\mtU(\co)$. To keep the notation compact, by the $p$-adic valuation of a (truncated) series in $\co$, we mean the minimum of the valuations of the coefficients. Similarly, $p$-adic valuation of a matrix refers to the minimum of the valuations of its entries, so that for example
\begin{align}
\begin{split}
	\ord(f_i^{[Cp]}(\varphi)) &\; \defineas \; \min \big\{ \ord(c_{i,n}) \; \big| \; n \in \{0,\dots, \left \lceil{Cp} \right \rceil \} \big\}~,\\ 
	\ord(\mtE(\co)) &\; \defineas \; \min \big\{\ord([\mtE(\co)]^{i,k}) \; \big| \; i,k \in \{0,\dots,b-1\} \big\}~.
\end{split}
\end{align}
To derive a bound for these valuations, it is useful to write the recurrence \eqref{eq:period_recurrences} in the form
\begin{align} \label{eq:period_recurrence_refined}
c_{i,n} \= \sum_{j = 0}^i \sum_{k = 1}^{N} \frac{P_{j,n-k}(n)}{n^{b+i}} \; c_{j,n-k}~, 
\end{align}
where $P_{i,k}(x) \in \IZ[x]$, and we have applied \eqref{eq:period_recurrences} to express the terms on the right-hand side in terms of the period coefficients $c_{j,m}$ with $j<i$ and $m<n$.

Taking the $p$-adic valuation of the recurrence \eqref{eq:period_recurrence_refined}, we obtain a simple recursive bound:
\begin{align} \label{eq:period_valuation_recurrence}
\ord(c_{i,n}) \, \geq \, \min \big\{-(b+i) \, \ord(n) + \ord(c_{j,n-k}) \; \big| \; j = 0,\dots,i,~ k = 1,\dots, N \big\}~, 
\end{align}
where we have used that $\ord(P_{j,m}(n)) \geq 0$.

Recalling that the $\co$-expansions can be truncated so that $n \leq \left \lceil{C p}\right \rceil$ we can estimate the $p$-adic valuation of $c_{i,n}$ from the above recurrence as
\begin{align} \label{eq:recurrence_sum}
\ord(c_{i,n}) \; \geq \; - (b+i) \, \sum_{m=0}^{\left \lceil{Cp}\right \rceil } \ord(m)~, \qquad \text{for} \qquad n \leq \left \lceil{C p}\right \rceil~,
\end{align}
Moreover, for primes $p > C$, this sum contains only terms with $m<p^2$, and hence 
\begin{equation}
    \ord(c_{i,n}) \; \geq \; - (b+i)\left\lceil C \right\rceil
\end{equation}
from which it follows that
\begin{align}
\ord([\mtE(\co)]^{i,k}) \geq - (b+i)\left\lceil C \right\rceil~.
\end{align}
In the expression \eqref{eq:U_matrix_definition} for the matrix $\mtU_p(\co)$, $\mtE(\co)$ appears multiplied by $\mtU_p(0)$. Due to the presence of the constant diagonal matrix $\diag(1,\dots,p^{b-1})$ in the definition of $\mtU_p(0)$, eq.~  \eqref{eq:U(0)}, the $i$'th row of the matrix $\mtE(\co)$ gets multiplied by $p^i$. Therefore
\begin{align}
\ord \left(\mtU_p(0) \, \mtE(\co) \right) \; \geq \; - (2b-1)\left\lceil C \right\rceil + b-1 + \min \{ \; \ord(\alpha_i) \; \mid \; i = 0,\dots,b-1\}~,
\end{align}
where $\alpha_i$ are the coefficients appearing in the matrix $\mtU_p(0)$. To obtain a similar expression for the inverse matrix $\mtE(\co^p)^{-1}$, truncated to order $\lceil Cp \rceil$ in $\varphi$, we note that the series expansions in $\mtE(\co^p)^{-1}$ involve only coefficients $c_{i,n}$ with $n \leq \left \lceil{C}\right \rceil$. For primes $p> \left \lceil{C}\right \rceil$, it follows from eq.~\eqref{eq:period_valuation_recurrence} that for these coefficients $\ord(c_{i,n}) \geq 0$. Therefore, recalling the inversion formula eq.~\eqref{eq:E_inversion} for $\mtE(\co)$,
\begin{align}
\ord(\mtE(\co^p)^{-1}) \= \ord(\mtW(\co^p)^{-1})~. 
\end{align}
Since $\mtW(\co)$ is a matrix of prime-independent rational functions, we can expand the entries of $\mtW(\co)^{-1}$ as series in $\varphi$ to order $\lceil C \rceil$, and the valuation of this matrix gives $\ord(\mtW(\co^p)^{-1})$, truncated to order $\lceil Cp \rceil$ in $\varphi$. In particular, this valuation vanishes for all but finitely many primes. 

Putting everything together, we conclude that to compute $\mtU_p(\co)$ to order $p^B$, it suffices to compute the inverse matrix $\mtE(\co^p)^{-1}$ to $p$-adic order
\begin{align}
B + (2b-1) \left \lceil{C}\right \rceil - b + 1 - \min \{ \; \ord(\alpha_i) \; \mid \; i = 0,\dots,b-1\}~.
\end{align}
This guarantees that, after multiplication by the entries of $\mtE(\co)$, all contributions of valuation at most $B$, are included. Similarly, the series in $\mtU_p(0)\mtE(\co)$ need to be computed to $p$-adic order
\begin{align}
B - \ord(\mtW(\co^p)^{-1})~.
\end{align}

\subsection{\texorpdfstring{$p$}{p}-adically truncated recurrence} \label{sect:truncated_recurrence}

To obtain the greatest benefit from using the $p$-adically truncated coefficients, it is not enough to compute the rational coefficients and then truncate them, as this would still require computing large rational numbers, and hence the memory usage would not be improved, making it impossible to study larger primes using this method. We therefore look for a recurrence which allows us to compute these coefficients directly.

Fix a prime $p$ and an \textit{initial} $p$-adic accuracy $A$. Then we define the \textit{$p$-adically truncated recurrence relation}
\begin{align} \label{eq:truncated_recurrences}
c_{i,n}^{(A)} \; \defineas \; \sum_{j = 0}^i \sum_{k = 1}^{N} \frac{P_{j,n-k}(n)}{n^{b+i}} \; c_{j,n-k}^{(A)}~ \mod p^A~,
\end{align}
which has the same form as the recurrence \eqref{eq:period_recurrence_refined} satisfied by the exact rational period coefficients but with the difference that after every step, the $p$-adic expansion of the resulting coefficient is truncated to order $A$. The initial conditions are also chosen to be equal to those of the recurrence \eqref{eq:period_recurrence_initial_conds} for the exact coefficients $c_{i,n}$: 
\begin{align}
c_{i,0}^{(A)} \= \delta_{i,0}~, \qquad c_{i,n}^{(A)} \= 0 \text{ for } n<0~. 
\end{align}
We call the solutions $c_{i,n}^{(A)}$ to this recurrence the \textit{approximate period coefficients} to differentiate them from the exact rational period coefficients $c_{i,n}$ modulo $p^A$. We define the approximate period series $f_i^{(A)}(\co)$ and period matrices $\mtE^{(A)}(\co)$ by formulae \eqref{eq:period_series_definition} and \eqref{eq:E_tilde_definition} but with the period coefficients $c_{i,n}$ replaced by the approximate coefficients $c_{i,n}^{(A)}$. 

The utility of the coefficients $c_{i,n}^{(A)}$ is that these provide a $p$-adic approximation of the exact rational coefficients $c_{i,n}$, and, by definition, they never grow larger than $p^A$. Thus they are significantly faster and less memory-intensive to compute than the $c_{i,n}$. Note, however, that the coefficients $c_{i,n}^{(A)}$ do not in general agree with $c_{i,n}$ modulo $p^A$, since division by $n^{b+i}$ in the recurrence relation \eqref{eq:truncated_recurrences} can lead to loss of accuracy. Therefore, to be able to reliably compute with these, we need to find a lower bound for how accurate approximation the coefficients $c_{i,n}^{(A)}$ provide given an initial accuracy $A$. Define the accuracy of the coefficient $c_{i,n}^{(A)}$ as
\begin{align} \label{eq:accDef}
\acc\big(c_{i,n}^{(A)}\big) \; \defineas \; \ord\big(c_{i,n}-c_{i,n}^{(A)}\big)~.
\end{align}
The recurrence relation \eqref{eq:truncated_recurrences} can be utilized to derive a recurrent formula for a lower bound for these accuracies. It implies the congruence
\begin{align}
c_{i,n} - c_{i,n}^{(A)} \; \equiv \; \sum_{j = 0}^i \sum_{k = 1}^{N} \frac{P_{j,n-k}(n)}{n^{b+i}} \; \big(c_{j,n-k} - c_{j,n-k}^{(A)} \big) \mod p^A~,
\end{align}
so that we can estimate the accuracy of these coefficients by
\begin{align} \label{eq:accuracy_recurrence}
\acc\big(c_{i,n}^{(A)}\big) \; \geq \; \min \left( \left\{  \acc\big(c_{j,n-k}^{(A)}\big) -(b+i) \, \ord(n) \; \middle | \; j = 0,\dots,i~,~~k = 1,\dots,N \right\}, A \right)~,
\end{align}
where $A$ appears due to the truncation to $p$-adic accuracy $A$.

Following the same strategy as in the previous \Cref{sect:series_truncation}, this can be used to derive a simple universal formula for lower bounds for the accuracies of the matrices $\mtU_p(0)\mtE^{(A)}(\co)$ and $\mtE^{(A)}(\co^p)^{-1}$. The matrix $\mtE^{(A)}(\co)$ involves coefficients $c_{i,n}^{(A)}$ with $n \leq \left \lceil{Cp}\right \rceil$. For these, when $p>C$,
\begin{align}
\acc\big( \, c_{i,n}^{(A)} \, \big) \; \geq \; A - (b+i) \, \sum_{m=0}^{\left \lceil{Cp}\right \rceil } \ord(m) \; \geq \; A - (b+i) \left\lceil C \right\rceil~.
\end{align}
Taking into account, as in the previous section, the factors of $p$ appearing in $\mtU_p(0)$, we arrive at the accuracy
\begin{align}
\text{acc}\left(\mtU_p(0) \, \mtE^{(A)}(\co)\right) \; \geq \; A - (2b -1) \left\lceil C \right\rceil + b-1 + \min \{ \; \ord(\alpha_i) \; \mid \; i = 0,\dots,b-1\}~.
\end{align}
To obtain a similar bound for the accuracy of the matrix $\mtE^{(A)}(\co^p)^{-1}$, note that this matrix involves only coefficients $c_{i,n}^{(A)}$ with $n \leq \left \lceil{C}\right \rceil$. Therefore, for primes $p> \left \lceil{C}\right \rceil$, the valuations $\ord(n)$ appearing in the recursive formula \eqref{eq:accuracy_recurrence} vanish. Hence,
\begin{align}
	\acc\big(\mtE^{(A)}(\co^p)^{-1}\big) \geq A + \ord\big(\mtW(\co^p)^{-1}\big)~.
\end{align}
Comparing these to the requisite accuracies derived for the matrices $\mtE^{(A)}(\co^p)^{-1}$ and $\mtU_p(0)\mtE^{(A)}(\co)$ in the previous section, we see that for primes $p > \left \lceil{C}\right \rceil$, the requisite accuracies are attained for both matrices if we pick\footnote{\label{foot:bound} It is easily possible provide a sharper bound if we do not insist on it being prime or operator-independent. For instance, by noting that the sum appearing in eq.~\eqref{eq:recurrence_sum} can be written as as $\left \lfloor \lceil C p \rceil/p \right \rfloor$, one obtains a sharper bound
\begin{align} \nonumber
	A = \max(B + (2b-1)\left \lfloor \lceil C p \rceil/p \right \rfloor - b +1 - \ord(\mtW(\co^p)^{-1}) - \min \{ \; \ord(\alpha_i) \; \mid \; i = 0,\dots,b-1\}, B - \ord(\mtW(\co^p)^{-1}))~.
\end{align} 
We could also use the $N$-integrality property of the coefficients $c_{0,n}$ explained in \Cref{app:Integrality} to improve the bounds for primes $p \nmid N$. Some such obvious improvements are implemented in the attached code (see \Cref{app:PFLFunction}) but are not discussed here to avoid proliferation of detailed analysis of different cases.}
\begin{align} \label{eq:A_lower_bound}
	A \= B + (2b-1) \left \lceil{C}\right \rceil - b +1 - \ord\left(\mtW(\co^p)^{-1}\right) - \min \{ \; \ord(\alpha_i) \; \mid \; i = 0,\dots,b-1\}~.
\end{align}

\section{Examples and Applications}
\label{sec:examples}
\noindent In this section, we give a few brief examples of computations which are now made possible by the $p$-adic truncated recurrence, and indicate some applications where such computations are of practical use. These examples also work as non-trivial consistency checks of our method.
\subsection{Computing the Euler factors for the family of mirror quintic threefolds at large primes}
\label{subsec:pBigQuintic}
To give a non-trivial consistency check of our method, we compute the Euler factors for $p = 2^{20}-3 = 1048573$ associated to the Calabi--Yau operator with the AESZ database \cite{AESZ} label \texttt{4.1.1},
\begin{align}
\cL \= \theta^4 - 5 \co (5 \theta + 1) (5 \theta + 2) (5 \theta + 3) (5 \theta + 4)~,
\end{align}
which corresponds to the family of mirror quintic threefolds whose middle cohomologies have a Hodge structure of type $(1,1,1,1)$. From the Weil conjectures, we have that $|a_p^{(2)}|\leq 6 p^3 < p^4$, where the last inequality holds for $p \geq 7$ so that for these primes $B=4$. In \cite{CdlOvS21}, it was found by extensive numerical experimentation that for this operator $C = 4/5$. Therefore, by eq.~\eqref{eq:A_lower_bound}, we can take the initial accuracy for our computations as $A = 8$.\footnote{Using the bound of \Cref{foot:bound}, we could in fact take $A=4$. However, we use $A=8$ in this section to demonstrate that the universal bound can be used to obtain results for large $p$ within reasonable time and memory constraints.}

Using a single core of an \texttt{Apple M5 Chip} (4.61 GHz), we can compute the Euler factors in 22 minutes, with a memory footprint of approximately 8500 MB. For instance, we find that at $\co = -1$
\begin{align}
R^{(3)}_{ 2^{20}-3}(\cL,T) \= 1-1576492860 \, T + 2672053179370 \, p \, T^2 - 1576492860 \, p^3 \, T^3 + p^6 \, T^4~.
\end{align}
This particular Euler factor was computed in ref.~\cite{controlledreduction} using \texttt{controlledreduction}, and the result given in ref.~\cite{controlledreduction} agrees with the above.

Note, however, that the method used in ref. \cite{controlledreduction} is somewhat orthogonal to our method: the method of ref. \cite{controlledreduction} is concerned with hypersurfaces in toric varieties (such as the quintic threefold with $h^{2,1} = 101$) and is most useful for computing the Euler factor for a single value of the moduli $\varphi$. In contrast, our method is only practically applicable for manifolds with a low number of complex structure moduli (such as the mirror quintic threefold with $h^{2,1}=1$), but is not restricted to a particular geometric construction. Additionally, our method is at its most efficient when computing the Euler factors for the whole family. Indeed, in 22 minutes, we also computed the other $1048571$ Euler factors at $p = 2^{20}-3$, making the time used per Euler factor approximately 0.0013 seconds.

We can also use the mirror quintic example to compare the memory usage of two variants of the deformation method algorithm when used to compute the Euler factors $R^{(3)}_p(\cL,T)$ for all primes $p \leq p_{\text{max}}$ for some $p_{\text{max}}$: The first computes the series $f_i(\co)$ with their full rational coefficients $c_{i,n}$, and uses these for the computation of the matrix $\mtU_p(\co)$. The second version uses the full machinery discussed in this paper, finding, for each prime $p$, the series approximate period coefficients $c_{i,n}^{(8)}$ to $p$-adic accuracy $8$ by using the $p$-adically truncated recurrence of \Cref{sect:truncated_recurrence}. 

In \Cref{fig:memory}, we display the peak memory usage of the computation of the polynomials $R_p^{(3)}(X_\co,T)$ using rational periods and periods obtained from the $p$-adically truncated recurrence to serve as a benchmark. The memory usage of the computation done using the full rational periods grows non-linearly, which is in stark contrast to the $p$-adically truncated recurrence, where the memory usage grows approximately linearly. For instance, computing the Euler factor at prime $p=23321$ takes 515.37 MB of memory using the method presented in this paper, whereas the same computation with rational periods requires over 6.8 GB of memory.

Finally, we present comparison of the time taken to compute the Euler factors using the full rational periods versus the periods obtained from the $p$-adically truncated recurrences. We compute the Euler factors for all $p \leq p_\text{max}$ using both methods. The rational method benefits from the fact that the full rational periods can be re-used for every $p \leq p_\text{max}$, whereas the $p$-adic recurrences must be computed anew for every prime. Consequently, the rational method is slightly faster at small primes. As primes grow, the $p$-adic method becomes faster. We present a comparison of the timings in \Cref{fig:timing}. Note, however that the timings presented here only compare primes which are relatively small due to the memory being the bottleneck for computing the rational periods. In addition, one of the benefits of the $p$-adic method is its parallelizability, which can be used to significantly speed up the computation.

\begin{figure}[h!]
    \centering
    \includegraphics[width=\linewidth, height=7cm, keepaspectratio]{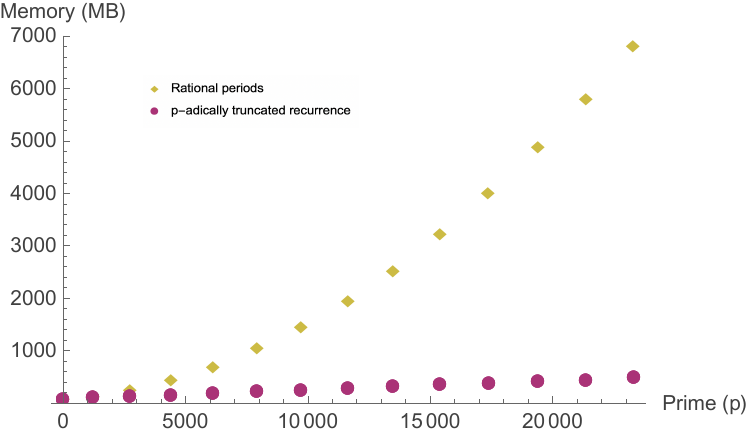}
    \capt{\linewidth}{fig:memory}{Comparing the memory usage for computing the Euler factors of the mirror quintic, using rational periods and $p$-adically truncated recurrence. For each prime $p$, the peak memory usage during the computation of $R_p^{(3)}(\cL,T)$ is displayed.}
\end{figure}

\begin{figure}[h!]
    \centering
    \includegraphics[width=\linewidth, height=7cm, keepaspectratio]{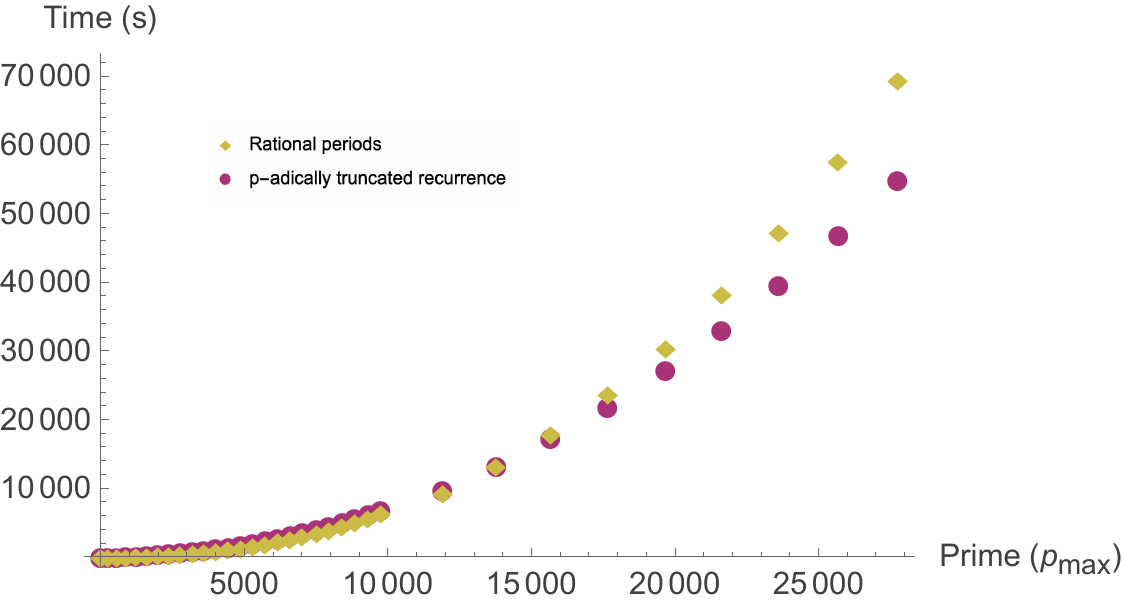}
    \capt{\linewidth}{fig:timing}{Comparing the timings of Euler factor computations using rational periods and periods obtained from the $p$-adically truncated recurrence. For each value $p_\text{max}$ displayed, we computed the Euler factors $R^{(3)}_p(\cL,T)$ for all $p\leq p_\text{max}$. Due to the memory usage being the bottleneck for the rational computation, we have only timed the computations for relatively small primes. It should also be emphasized that the timings of the $p$-adically truncated method can be substantially improved by running different primes in parallel; this speed-up is algorithmically infeasible for the rational computation.}
\end{figure}

\subsection{Using statistical properties of the Frobenius traces}
Being able to compute the Frobenius traces for large primes also allows us to study and make use of their statistical properties such as Sato--Tate distributions. The Sato--Tate conjecture, proven in \cite{SatoTateProofI, SatoTateProofII, SatoTateProofIII} (see for example \cite{Sutherland2019a} and references therein for an accessible introduction) shows that for elliptic curves, the complex multiplication property can be linked to distributions of the coefficients $a_p$: the conjecture states that for elliptic curves $E/\IQ$ without complex multiplication, the values $a_p/\sqrt{p}$, which satisfy $|a_p/\sqrt{p}| \leq 2$, are equidistributed with respect to the measure
\begin{align}
	\dd \mu \= \frac{1}{2\pi} \sqrt{4-t^2} \; \dd t~.
\end{align}
The equidistribution means that for any interval $[\alpha,\beta] \subseteq [-2,2]$, the following limit holds
\begin{align}
 \lim_{C \to \infty} \frac{\#\{p \leq C \; | \; a_p/\sqrt{p} \in [\alpha,\beta]\}}{\#\{p \leq C\}} \= \frac{1}{2\pi} \int_{\alpha}^\beta \sqrt{4-t^2} \; \dd t~.
\end{align}
For curves $E/\IQ$ with complex multiplication, the measure is given by
\begin{align} \label{eq:CM_measure}
 \dd \mu_\text{CM} \= \frac{1}{2 \pi} \frac{1}{\sqrt{4-t^2}} \; \dd t + \frac{1}{2} \delta(t)	\; \dd t~,
\end{align}
where $\delta(t)$ is the Dirac distribution.

Computing numerically distributions of Frobenius traces unlocks concrete applications, for instance, in mathematical physics. In ref.~\cite{CdlOEvS20}, the factorization of the Euler factors $R^{(3)}_p(\cL,T)$ over $\IQ$ into two quadratic factors was used to find \textit{rank-two attractor varieties}, Calabi--Yau threefolds for which the Hodge structure of the four-dimensional middle cohomology splits, over $\IQ$, to two two-dimensional subspaces, of Hodge types $(3,0)+(0,3)$ and $(2,1)+(1,2)$, respectively. For these varieties, the Galois representation $\rho: \Gal(\overline{\IQ}/\IQ) \to \GL(H^3(X_\co, \IQ_\ell))$ (conjecturally) splits as $\rho = \rho_4 \oplus \rho_2$, where $\rho_4$ is weight-4 modular and $\rho_2$ weight-2 modular. Such points are of significant interest in string theory compactifications owing to their relation to black holes and flux vacua (see for example refs.~\cite{KlemmEtAl,CdlOEvS20,Candelas:2023yrg,JKK23,Grimm:2024fip,Blesse:2025dch,Ducker:2025wfl,Kachru:2020abh,Kachru:2020sio}). Computing trace distributions not only gives us an alternative way of finding rank-two attractor varieties, but also allows us to probe finer properties of the Galois representations, and thus of Hodge structures, such as the complex multiplication property, which cannot be found by studying factorizations over $\IQ$ alone.\footnote{In this case, one could also study factorizations over field extensions. However, unless one knows a priori over which field extension the factorizations should be expected to occur, this gets quickly cumbersome due to having to take into account the ramification structure of primes over different field extensions.}

To give a toy example of using trace distributions in practice to find interesting varieties, let us consider the question of finding singular (also termed \textit{attractive} in mathematical physics contexts \cite{Moore}) K3 surfaces among a one-parameter family of K3 surfaces.\footnote{Of more interest to physicists is the question of finding analogous Calabi--Yau fourfolds due to their relation to M-theory vacua. However, discussing this more involved case would take us outside of the scope of the present note, so we aim to return to this topic in a future study~\cite{Kuusela:2026mst}.} As our example, we take the one-parameter family of K3 surfaces $S_\co$ studied in ref.~\cite{Verrill} which are defined as the toric compactifications of vanishing loci of
\begin{equation}
    (X_0 + X_1 + X_2 + X_3) \left( \frac{1}{X_0} + \frac{1}{X_1} + \frac{1}{X_2} + \frac{1}{X_3} \right) \= \frac{1}{\co} \qquad \mathbb{P}^3 \, \setminus \, \{ X_{i} \=  0 \}~.
\end{equation}  
It has been shown in ref.~\cite{Verrill} that the Picard number of a generic member of this family is $19$. Therefore, the Picard--Fuchs operator is of order 3 and reads
\begin{align} \label{eq:HV_PF_operator}
\cL \= \theta^3 + 64 \, \co^2 \, ( \theta +1 )^3 - 2 \, \co \, ( 2 \theta + 1 ) ( 5 \theta \, ( \theta +1  ) +2 )~.
\end{align}
Thus, the corresponding Euler factors are of the form
\begin{align} \label{eq:K3_Euler_factor_factorised}
R^{(2)}_p(\cL,T) \= 1 + a^{(1)}_p(\co) T + a^{(2)}_p(\co) T^2 + a^{(1)}_p(\co)/a^{(2)}_p(\co) p^4 T^3 \= (1 \pm p T) (1 + c_p(\co) T + p^2 T^2)~. 
\end{align}

Since this polynomial factorizes always over $\IQ$, one cannot use the analogues of the techniques of ref.~ \cite{CdlOEvS20} to find singular K3 surfaces. However, we can make progress by using the fact that singular K3 surfaces defined over $\IQ$ have a two-dimensional Galois representation acting on the transcendental lattice, so that the Galois representation $\rho: \Gal(\overline{\IQ}/\IQ) \to \GL(H^2(S_\co,\IQ_\ell))$ splits into 1-dimensional representation $\rho_1$ and a two-dimensional representation $\rho_2$, $\rho = \rho_1 \oplus \rho_2$. Furthermore, the representation $\rho_2$ is modular and of complex multiplication type in the sense that there exists a Hecke eigenform of weight 3 with complex multiplication such that the Hecke eigenvalues agree with the traces of $\rho_2$ \cite{Livne1995a, Elkies2013a}. The traces of such a representation $\rho_2$ are equidistributed with respect to the measure $\dd \mu_{\, \text{CM}}$ defined in eq.~\eqref{eq:CM_measure}. If the representation $\rho_1$ is trivial, its traces are trivially equidistributed with respect to the Dirac measure $\delta(x+1) \dd x$. Then the distribution $C(x)$ of traces of $\rho$ can be found by computing the product
\begin{align} \label{eq:semicircle_measure}
\dd \mu_{\, \text{CM}} \, \delta(x+1) \dd x \= \frac{1}{2\pi} \frac{1}{\sqrt{4-(x+1)^2}} \dd x + \frac{1}{2} \, \delta(x+1) \dd x \asdefine C(x) \dd x~.
\end{align}
Similarly, the non-trivial representation $\rho_1$ has the corresponding distribution $\frac{1}{2}(\delta(x-1) + \delta(x+1))$, and the product with the distribution $ \mu_{\, \text{CM}}$ gives the distribution $F(x)$ which was called the ``flying Batman'' distribution in ref.~\cite{Chen2024a}\footnote{In contrast to ref.~\cite{Chen2024a}, we have included the Dirac distributions in the definition.}
\begin{align}  \label{eq:flying_Batman_distribution}
F(x) \= \begin{cases}
\frac{1}{4 \pi \sqrt{3 - 2x - x^2}} + \frac{1}{4 \pi \sqrt{3 + 2x - x^2}} & |x| < 1~,\\
\frac{1}{4 \pi \sqrt{3 + 2|x| - x^2}} + \frac{1}{4}(\delta(x-1) + \delta(x+1)) & 1 \leq |x| \leq 3~,\\
0 & \text{otherwise.}
\end{cases}
\end{align}
In the case the surface is not of CM-type, we expect to find either the pushforward of the Haar measure on $\text{O}(3)$, which is the ``Batman'' distribution $B(x)$ \cite{Saad2023a,Ono2023a}, or that of the Haar measure on $\SO(3)$ which, to continue the Batman metaphor, we call here the ``wing'' distribution $W(x)$:
\begin{align} \label{eq:non-CM_distributions}
B(x) \= \begin{cases}
\frac{3+x}{4 \pi \sqrt{3 - 2x - x^2}} + \frac{ 3-x}{4 \pi \sqrt{3 + 2x - x^2}} & |x| < 1~,\\
\frac{3-|x|}{4 \pi \sqrt{3 + 2|x| - x^2}} & 1 \leq |x| \leq 3~,\\
0 & \text{otherwise.}
\end{cases} \qquad W(x) \= \frac{1}{2 \pi} \sqrt{\frac{3+x}{1-x}}~.
\end{align}
To be able to numerically classify the distributions of Frobenius traces for the first $10000$ primes we have computed, we can compute the first few moments of the distributions:
\begin{align}
M(\co) \defineas \frac{1}{10000} \sum_{p}(1, a_p^{(1)}(\co), a_p^{(1)}(\co)^2, \dots, a_p^{(1)}(\co)^5)~.
\end{align}
We can compare these to the corresponding moments of the above distributions, which are given by
\begin{align}
\begin{split}
M_B &\= (1,0,1,0,3,0)~, \qquad M_W \= (1,0,1,-1,3,-6)~,\\
M_F &\= (1,0,2,0,10,0)~,  \qquad M_{C} \= (1,-1,2,-4,10,-26)~.
\end{split}
\end{align}
To determine which (if any) of the expected distributions in eqs.~\eqref{eq:semicircle_measure}, \eqref{eq:flying_Batman_distribution}, \eqref{eq:non-CM_distributions} the distribution of Frobenius traces at $\co$ resembles, we can compute the squared distance $(M(\co) - M_D)^2$ for each distribution ${D \in \{B,W,F,C\}}$. We compute these for all values $\varphi = r/s$ with $|r|,|s| \leq 250$. We find that every distribution $M(\co)$ has squared distance of less than 2 to one of the distributions  $D \in \{B,W,F,C\}$: 75 937 of these correspond to the ``Batman'' distribution $B$ in eq.~\eqref{eq:non-CM_distributions}, and 151 to the ``wing'' distribution $W$ in eq~\eqref{eq:non-CM_distributions}. For $\co = 1, -1/2, -1/8, -1/32, 1/64$, we find the ``flying Batman'' distribution \eqref{eq:flying_Batman_distribution}, and for $\co = 1/8$ the semicircle distribution \eqref{eq:semicircle_measure}. We expect these six cases to correspond to K3 surfaces of complex multiplication type. Indeed, in all of these cases, we can identify the coefficients $c_p(\co)$ in the Euler factor $R_p^{(2)}(\cL,T)$ of eq.~\eqref{eq:K3_Euler_factor_factorised} with Hecke eigenvalues of weight-3 modular forms, at least to the first 7000 eigenvalues. We have listed the explicit identifications in \Cref{tab:CM_modular_forms}. In \Cref{fig:CM_HVK3,fig:Non_CM_HVK3} four histograms of Frobenius traces we have computed are displayed, one corresponding to each of the four possible distributions discussed above.

\begin{table}[h!]
\begin{tabular}{|l|c|c|c|c|c|c|}
\hline
$\co$    & $1$ & $-1/2$ & $1/8$ & $-1/8$ & $-1/32$ & $1/64$  \\ \hline
label    & \texttt{15.3.d.b} & \texttt{12.3.c.a} & \texttt{8.3.d.a} & \texttt{24.3.h.b} & \texttt{12.3.c.a} & \texttt{15.3.d.b} \\ \hline
\end{tabular}
    \vskip7pt
    \capt{6.5in}{tab:CM_modular_forms}{The LMFDB labels of modular forms associated to L-functions.}
    \vskip-20pt	
\end{table}

\begin{figure}[h!]
	\centering
	\begin{center}
		\includegraphics[width=\linewidth, height=7cm, keepaspectratio]{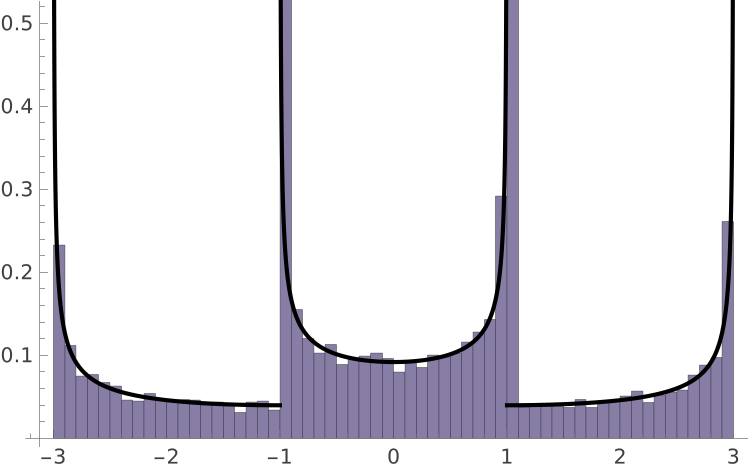}
        \vskip30pt
        \includegraphics[width=\linewidth, height=7cm, keepaspectratio]{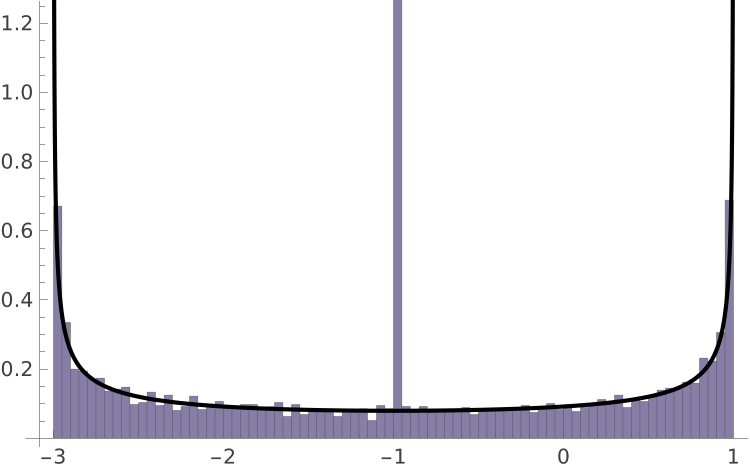}
	\end{center}
        \place{3.2}{0.15}{$a_p^{(1)}/p$}
        \place{0.8}{2.2}{\rotatebox{90}{Relative abundance}}
        \place{3.2}{3.3}{$a_p^{(1)}/p$}
        \place{0.8}{5.5}{\rotatebox{90}{Relative abundance}}
        \capt{\linewidth}{fig:CM_HVK3}{The normalized histogram displaying the distribution of the normalized Frobenius traces $a_p^{(1)}/p$ for the order 3 degree 2 Calabi--Yau operator \eqref{eq:HV_PF_operator} at $\varphi = 1$ (upper figure) and $\varphi = 1/8$ (lower figure). We include the traces corresponding to the first 10000 primes in 70 bins. We see that the figures follow approximately the ``flying Batman'' \eqref{eq:flying_Batman_distribution} and shifted semicircle \eqref{eq:semicircle_measure} distributions, respectively, which we have overlaid on the histogram. These distributions are expected when the corresponding K3 surface has complex multiplication. Note that the figures do not display the bars fully at $\pm 1$, which together account for approximately half of the traces of the respective figures.}	
\end{figure}

\begin{figure}[h!]
	\centering
	\begin{center}
		\includegraphics[width=\linewidth, height=7cm, keepaspectratio]{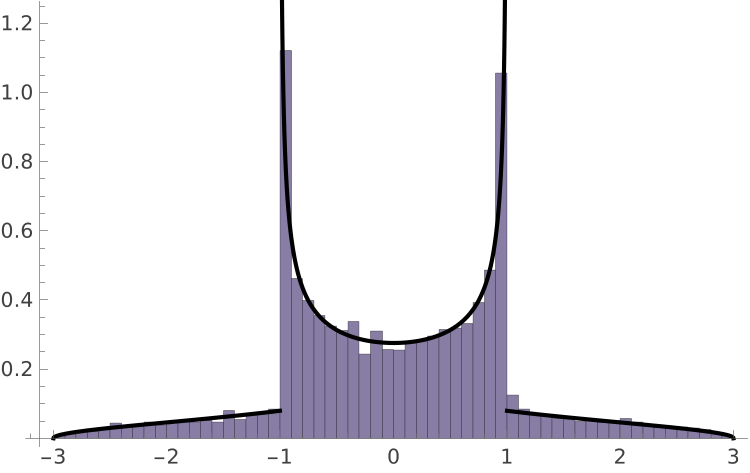}
        \vskip30pt
        \includegraphics[width=\linewidth, height=7cm, keepaspectratio]{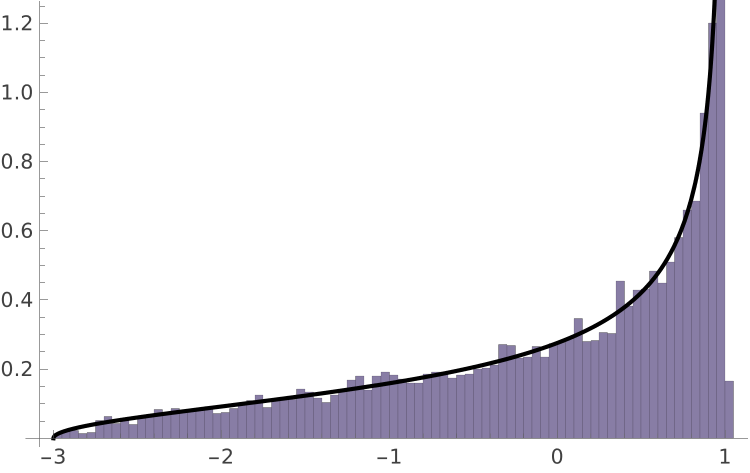}
        \vskip30pt
        \capt{\linewidth}{fig:Non_CM_HVK3}{The normalized histogram displaying the distribution of the normalized Frobenius traces $a_p^{(1)}/p$ for the order 3 degree 2 Calabi--Yau operator \eqref{eq:HV_PF_operator} at $\varphi = 3$ (upper figure) and $\varphi = 2/7$ (lower figure). We include the traces corresponding to the first 10000 primes in 70 bins. We see that the figures follow approximately the ``Batman'' and ``wing'' distributions of eq.~\eqref{eq:non-CM_distributions}, respectively, which we have overlaid on the histogram. These distributions are the pushforwards of $\text{O}(3)$ and $\SO(3)$ Haar measures, and correspond to cases where the K3 surface does not have complex multiplication.}	
	\end{center}
        \place{3.2}{1.4}{$a_p^{(1)}/p$}
        \place{0.8}{3.5}{\rotatebox{90}{Relative abundance}}
        \place{3.2}{4.6}{$a_p^{(1)}/p$}
        \place{0.8}{6.7}{\rotatebox{90}{Relative abundance}}
\end{figure}

\subsection{Predicting Hecke eigenvalues of paramodular forms}
It is expected that the Calabi--Yau motives of Hodge type $(1,1,1,1)$ correspond to Siegel paramodular forms of genus $g = 2$, weight $3$, and some level $N$. These modular forms admit Hecke theory, and the Hecke eigenforms are characterized by the eigenvalues $\lambda_{p,1}$ and $\lambda_{p,2}$ of two Hecke operators $T_{p,1}$ and $T_{p,2}$ for every good prime $p$. For levels $N < 1000$ and primes $p < 200$, these eigenvalues have been tabulated in ref.~\cite{omf5data} using the relation to orthogonal modular forms \cite{Eran2024a,Rama2020a,Dummigan2024a} (see also ref.~\cite{qomf} for a database of the eigenvalues $p < 1000$ for levels $N < 100$). Many examples of such Calabi--Yau motives and Siegel paramodular forms have been found experimentally using the deformation method already in ref.~\cite{GvS24}. 

The utility of the $p$-adic truncation in this process is that it allows us to easily compute primes significantly higher than $p < 1000$ included in the databases. Therefore, if we can find a Calabi--Yau type differential operator $\cL$ and modulus $\co$ such that the associated motive (conjecturally) corresponds to one of these Siegel paramodular forms, we can use the zeta function data to predict further eigenvalues. Specifically, the relation of the coefficients of $R_p^{(3)}(\cL,T)$ to the eigenvalues $\lambda_{p,1}$ and $\lambda_{p,2}$, in the conventions of ref.~\cite{Rama2020a},
\begin{align}
\lambda_{p,1} \= -a_p^{(1)}~, \qquad \lambda_{p,2} \= \frac{a^{(2)}_p-p-p^3}{p}~.
\end{align}
The relation of the Euler factors to $L$-functions give us strong consistency checks on the method. The associated $L$-function $L(s)$ is conjectured to satisfy a functional equation: if we denote the completed $L$-function by
\begin{equation}  \label{eq:completed_L-function}
    \Lambda(s) \= \left(\frac{N}{\pi^4}\right)^{s/2}\Gamma\left(\frac{s-1}{2}\right)\Gamma^2\left(\frac{s}{2}\right)\Gamma\left(\frac{s+1}{2}\right)L(\cL,s)~
\end{equation}
it should satisfy the functional equation
\begin{equation} \label{eq:L_functional_equation}
    \Lambda(s) \= \epsilon\Lambda(4-s)~.
\end{equation}
Using the Euler factors we have computed, we can form the approximate $L$-function
\begin{align} \label{eq:approximate_L-function}
    L_{p_{\text{max}}}(\cL,s) \; \defineas \; \prod_{p=2}^{p_\text{max}}R_p(\cL,p^{-s})^{-1}~,
\end{align}
where the product is over all primes between $2$ and $p_{\text{max}}$. 

Then, by computing the analytic continuation of this function, we can find the accuracy to which the functional equation \eqref{eq:L_functional_equation} is satisfied. This is done by using the function \texttt{checkfeq} in PARI/GP, which implements Dokchitser's method~\cite{dokchitser}. In particular, we compute the binary accuracy to which the approximate $L$-function satisfies the functional equation \eqref{eq:L_functional_equation}. This is defined as the least integer $\eta_0$ such~that
\begin{align}
    |\Lambda_{p_{\max}}(s)-\epsilon\Lambda_{p_{\max}}(4-s)| \leq 2^{\eta_0(p_{\text{max}})}~,
\end{align}
where $\Lambda_{p_{\max}}(s)$ is formed out of the approximate $L$-function $L_{p_{\text{max}}}(\cL,s)$ analogously to eq.~\eqref{eq:completed_L-function}.

For instance, we can consider the Calabi--Yau operator with the AESZ \cite{AESZ} database \texttt{newnumber} \texttt{4.2.5} \cite{CYcluster}:
\begin{align}
\cL \= \theta^4-2^{2} \co (2\theta+1)^2(11\theta^2+11\theta+3)-2^{4} \co^{2}(2\theta+1)^2(2\theta+3)^2~.
\end{align}
This operator has discriminant
\begin{align}
\Delta(\co) \= 256 \co^2 + 176 \co - 1~,
\end{align}
so for $\co = -1$, the discriminant is $79$. We therefore expect that the corresponding modular form is of level $N$ such that $79 \vert N$. Indeed, by comparing to the database \cite{qomf}, we see that the Euler factors we compute from this operator correspond exactly to those of the form with the label \texttt{79a}, as observed already in ref.~\cite{GvS24}. However, we can easily compute further Euler factors that give a prediction for further Hecke eigenvalues. We compute the first 10000 Euler factors, and by using the PARI/GP function \texttt{checkfeq}, we find that $\eta_0(104729) =-1889$.

To get an indication of whether the Euler factors are correct even at large primes, we can study the accuracy as a function of the maximum prime $p_\text{max}$ included in the product \eqref{eq:approximate_L-function}. Doing this, we indeed find that the accuracy increases as more Euler factors are included, which is consistent with these factors being computed correctly. We display the accuracy as a function of $p_{\text{max}}$ in \Cref{fig:L_functional_equation}.

\begin{figure}[H]
	\centering
	\begin{center}
		\includegraphics[width=9cm, height=4.5cm]{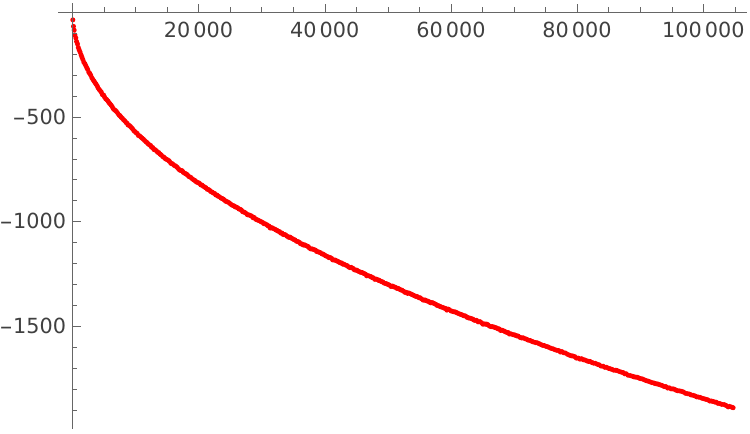}
	\end{center}
    \place{3.2}{2.0}{$p_{\max}$}
    \place{1.3}{1.0}{$\eta_0$}
	\capt{7in}{fig:L_functional_equation}{The binary accuracy $\eta_0$ to which the approximate $L$-function \eqref{eq:approximate_L-function} satisfies the functional equation \eqref{eq:L_functional_equation}. The horizontal axis displays the maximal prime $p_{\max}$ included in the approximate $L$-function, and the vertical axis the corresponding accuracy. We display the data for every 20th prime. Note that the accuracy increases as the number of Euler factors is increased. This strongly indicates that the Euler factors computed using our method are correct, as the accuracy would be expected to plateau after an incorrect Euler factor.}	
\end{figure}

\section*{Acknowledgements}
\noindent We thank John Voight for insightful discussions and comments which inspired writing this note. We are grateful to Philip Candelas and Xenia de la Ossa for discussions and comments on the manuscript. We wish to express our thanks for various interesting and useful conversations to Paul Blesse, Janis Dücker, Hans Jockers, Abhiram Kidambi, Albrecht Klemm, Sören Kotlewski, Eric Pichon-Pharabod, Gustavo Rama, Duco van Straten, Gonzalo Tornaría, and Wadim Zudilin. We thank Hans Jockers for providing computer resources which were used to compute much of the data discussed here. MMR and PK are supported by the Cluster of Excellence PRISMA+ (Project ID 390831469). MMR acknowledges the Sociedad Matemática Mexicana (SMM) together with the Fundación Sofía Kovalévskaia (SK) for their support through the Apoyo Sofía Kovalévskaia 2025, which contributed to the completion of this work. PK wishes to acknowledge support from Emil Aaltosen säätiö. PK is supported by the Leverhulme Trust grant ``Quantum geometry and arithmetic'' (NIP-2023-007). PK thanks the LPTHE of Sorbonne Université and CNRS, and in particular Boris Pioline, for the hospitality during the final stages of this work.
ML would like to thank Per Berglund for his support at the beginning of this project and Reiko Toriumi for her support at the end of this project.
The work of MS was supported in part by NSF grant PHY-2309456. MS would like to thank Liam McAllister for his support in pursuing this line of work, and Richard Nally for introducing him to this intersection between number theory and physics. We wish to express our gratitude for MITP for hospitality during the workshop ``The Arithmetic of Calabi--Yau Manifolds''. We would like to acknowledge the ICTP meetings ``Workshop on Number Theory and Physics'' and ``Workshop on Explicit Arithmetic Geometry'', during which much progress on the present work was made. 

\newpage
\appendix
\section{Integrality Conditions of Calabi--Yau Type operators} \label{app:Integrality}
\noindent In addition to the conditions listed in \Cref{sec:review}, differential operators of Calabi--Yau type satisfy the following \textit{integrality conditions}:
\begin{enumerate}[label=D\arabic*]
    \item The holomorphic solution of $\cL$ at $\co = 0$ is \textit{$N$-integral}, meaning that there exists an $N \in \mathbb{Z}$ such that $N^{m} \, c_{0,m} \in \mathbb{Z}$ for all $m \geq 0$. \label{cond:fundamental_period_N-integrality}
    \item All structure series $\beta_i$, $0 \leq i \leq b-1$, of $\cL$ at $\co = 0$ are $N$-integral. The \textit{i'th structure series} $\beta_i$ of $\cL$ is defined as
    \begin{equation} \label{eq:structure_series}
    \beta_{i} \; \defineas \; \mathcal{N}_{i}(\varpi^{i})^{-1} \, \in \, \mathbb{Q}\llbracket\co\rrbracket^{*}~,
    \end{equation}
    where $\mathcal{N}_{k} \in \mathbb{Q}\llbracket\co\rrbracket[\theta]$ for all $0 \leq k \leq b$, is a differential operator defined through the recurrence relation
    \begin{equation}
        \mathcal{N}_{0} \; \defineas \; 1 , \quad \mathcal{N}_{k+1} \; \defineas \; \theta  \, \frac{1}{\mathcal{N}_{k}(\varpi^{k})} \, \mathcal{N}_{k} ~,
    \end{equation}
    and satisfies $\mathcal{N}_{k}(\varpi^{i}) = 0 $ for all $0 \leq i \leq k-1 $~.
    \item The \textit{$q$-coordinate of $\cL$ at $\co = 0$} is $N$-integral. The $q$-coordinate, also known as the \textit{special coordinate}, is defined as\footnote{The $q$-coordinate can also be defined as the unique solution $q \in \co \, \mathbb{Q}\llbracket\co\rrbracket^{*}$ of the differential equation $\theta q = q \, \beta_{1}^{-1}$ with the initial condition $\partial_{\co}q(0) = 1$.}
    \begin{equation}
        q \= \exp\left( \frac{\varpi^{1}}{\varpi^{0}} \right) \, \in \, \co \, \mathbb{Q}\llbracket\co\rrbracket^{*}~,
    \end{equation}
\end{enumerate}

\section{Rational Structure Associated to Fourth-Order Calabi--Yau Operators}
\label{app:Rationality}
\noindent In this appendix, we give a few details of the construction of the rational structure associated to Calabi--Yau operators of order 4. For such operators which arise from a Calabi--Yau motive, the $\Gamma$-class conjecture (see e.g. refs.~\cite{Hori:2013ika,Iritani:2007a,Iritani2009a,Libgober:1999a,Halverson:2013qca}) implies that the period vector $\varpi = (\varpi^0,\dots,\varpi^3)$ is related to a period vector $\Pi$ with integral monodromies by
\begin{align}
\Pi \= \wh \rho \nu^{-1} \varpi ~,
\end{align}
where 
\begin{align} \label{eq:rho_Gamma_class}
\nu \; \defineas \; \text{diag}(1,2\pi\ii,(2\pi\ii)^2,(2\pi\ii)^3)~, \qquad \wh \rho \; \defineas \; \begin{pmatrix}
\chi \smallfrac{\zeta(3)}{(2\pi \ii)^3} & \smallfrac{c_2 \cdot \me}{24} & 0 & \me^3 \\
\smallfrac{c_2 \cdot \me}{24} & \smallfrac{\sigma}{2} & -\me^3 & 0\\
1 & 0 & 0 & 0\\
0 & 1 & 0 & 0\\
\end{pmatrix},
\end{align}
$\chi$ is the Euler characteristic of the associated mirror Calabi--Yau manifold, $c_2 \cdot \me$ its second Chern class, $\me^3$ the triple intersection number, and $\sigma \in \{0,1\}$. 

There are, however, known examples of Calabi--Yau type operators for which the above formula needs to be generalized (e.g. \cite{Katz:2022lyl,Katz:2023zan}). In ref.~\cite{PichonPharabod2025a}, a generalized formula for $\rho$ is proposed, which in our conventions can be written as:
\begin{align}
\rho \=
\begin{pmatrix}
\alpha  M  -\delta & L M & \alpha  S & S \\
 L & N  \frac{\sigma}{2}  & -S & 0 \\
 1 & 0 & 0 & 0 \\
 \frac{\alpha  N}{M} & N & 0 & 0 \\
\end{pmatrix}
\begin{pmatrix}
1 & 0 & 0 & 0 \\
 0 & 1 & 0 & 0 \\
 0 & 0 & 1 & 0 \\
 \smallfrac{\zeta(3)}{(2\pi \ii)^3} K & 0 & 0 & 1 \\
\end{pmatrix},
\end{align}
with $\alpha,\beta,\sigma,K,L,M,N \in \IQ$. Importantly for us, the multiplication by the first matrix with rational entries does not change the rational structure. Therefore, $\wh \rho$ and $\rho$ give rise to the same rational structure which determined by a single rational number $K$. In the case the $\Gamma$-class formula \eqref{eq:rho_Gamma_class} holds, this is given as $K = \chi/\me^3$. Following \cite{ZetaMultThreeFolds}, we conjecture that this change of basis gives a rational structure for all Calabi--Yau type operators of order 4 which have a rational structure. Further, we conjecture that the rational number $K$ appearing in this relation is equal to the coefficient $K$ in the matrix $\mtU_p(0)$. We are not aware of counterexamples to this conjecture. 

\section{\texttt{PFLFunction}: a \texttt{Sage}-compatible \texttt{Python} Package for Computing Euler Factors Associated to Calabi--Yau Type Operators}
\label{app:PFLFunction}
\noindent In this appendix, we present documentation for our \texttt{Python} package \texttt{PFLFunction}, which contains implementations of many algorithms described in the paper. The aim of the package is to make these as user-friendly as possible for one-parameter operators of Calabi--Yau type, with multiparameter functionality to appear in a later release.

\subsection{Downloading and installing}
\vskip-10pt
The package can be downloaded from
\begin{center}\href{https://github.com/PyryKuusela/PFLFunction/tree/main}{\texttt{https://github.com/PyryKuusela/PFLFunction}}
\end{center}

\noindent \textbf{Step 1:} In a terminal window, type:
\begin{lstlisting}[style=notebook]
    git clone https://github.com/PyryKuusela/PFLFunction.git
\end{lstlisting}
\noindent to download the package contents.

\noindent \textbf{Step 2:} Next, in the terminal window, change the current working directory via:
\begin{lstlisting}[style=notebook]
    cd PFLFunction
\end{lstlisting}
\noindent to \texttt{PFLFunction}.

\noindent \textbf{Step 3:} Next, create a virtual environment to have all of the requisite packages installed. This can be done using conda by typing the following into the terminal window:
\begin{lstlisting}[style=notebook]
    conda env create -f environment.yml
\end{lstlisting}
\noindent Note that in particular, the periods code depends on \texttt{SymPy} and \texttt{NumPy}. Power series computations are based on \texttt{FLINT} 
polynomial multiplication.

\noindent \textbf{Step 4:} Activate the virtual environment by typing the following into the terminal window:
\begin{lstlisting}[style=notebook]
    conda activate pflfunction
\end{lstlisting}

\noindent \textbf{Step 5:} To install the package contents in the virtual environment, type:
\begin{lstlisting}[style=notebook]
    pip install .
\end{lstlisting}
The package is now ready to use.

To use the package, we have included the folder \texttt{scripts}, which contains both an example notebook \texttt{example.ipynb}, as well as bash scripts that parallelize the zeta function computations across primes. The example notebook shows both how to use the main function of the code, \texttt{L\_functions} (detailed below), as well as how to create input data that can be used to run the provided bash script \texttt{parallel\_zeta.sh}. 

\begin{funcDefnN}[label=main]{\texttt{L\_functions}}
    \texttt{L\_functions(L,C,p,conifold\_locus,apparent\_sing\_locus,other\_sing\_locus,}\\ \hfill\texttt{K,label,pacc\_init=0,nadd=0)}\\
    The main function, which computes the polynomial $R_p^{(b-1)}(\mathcal{L},T)$ using the following data.\\
    \textbf{Arguments}
    \begin{itemize}
        \item \texttt{L}: A Calabi--Yau type operator. Currently, this is only supported for operators of order $3$ or $4$. Future support for operators of order $2$ and $5$ is planned.
        \item \texttt{C}: A rational number $C$ appearing in eq.~\eqref{eq:trunc_order}, determining the order $\lceil Cp \rceil$ in $\co$ to which the power series expansions are computed.
        \item \texttt{p}: The prime for which the polynomial $R_p^{(b-1)}(\mathcal{L},T)$ will be computed.
        \item \texttt{conifold\_locus}: The conifold locus as a \texttt{SymPy} expression.
        \item \texttt{apparent\_sing\_locus}: The apparent singular locus as a \texttt{SymPy} expression.
        \item \texttt{other\_sing\_locus}: The other singular loci as \texttt{SymPy} expressions. 
        \item\texttt{K}: The rational number $K$ specifying the rational structure associated to \texttt{L}. See \Cref{app:Rationality} for more details.
        \item \texttt{label}: A string used to name the output and log files.
        \item \texttt{pacc\_init}: The initial $p$-adic accuracy $A$ used for the $p$-adically truncated recurrence of \Cref{sect:truncated_recurrence}. If \texttt{nadd} is nonzero, the user should input a value to guarantee truncation; otherwise, the lower bound from eq.~\ref{eq:A_lower_bound} is computed.
        \item \texttt{nadd}: The number of terms that should vanish at the end of a series for the purposes of checking that the series terminates. This is a consistency check beyond the default expansion order $M(\mathcal{L},p)$; i.e. the relevant expressions will be computed to order $M(\mathcal{L},p) + \texttt{nadd}$ in $\co$. If \texttt{nadd} is nonzero, then the option \texttt{pacc\_init} should be set appropriately.
    \end{itemize}
    \textbf{Output}
    \begin{itemize}
        \item \texttt{R polynomial}: A tuple $[[p,\co_*,[a_p^{(k)}(\co_*)]] \texttt{ for $\co_*$ in $\mathbb{F}_p^*$}]$ of the prime $p$, the point $\co_*$, and the coefficients $a_p^{(k)}(\co_*)$ in the $R_p^{(b-1)}(\cL,T)$ as defined in eq.~\eqref{eq:a_p_definition}. As these are only yet implemented for $b=3,4$, the relevant coefficients printed for each value of $\co_*$ are $[a_p^{(1)},a_p^{(2)}]$. If a given value of $\co_*$ is a C-point, then this is printed as $[p,\co,[a_p^{(k)}(\co)],\text{C}]$. If the point $\co_*$ is another type of singular point of \texttt{L} $\!\!\mod p$, the output is  $[p,\co,[a_p^{(k)}(\co)],0]$    
    \end{itemize}
\end{funcDefnN}

\subsection{Setup and options}
To set up and use the bash script \texttt{parallel\_zeta.sh},
in a terminal window with the directory set to the \texttt{scripts} folder, type:
\begin{lstlisting}[style=notebook]
    chmod +x parallel_zeta.sh
\end{lstlisting}
\noindent to turn the bash script into an executable.
The bash script takes in six inputs: \texttt{operator\_name}, \texttt{primes}, \texttt{scaling}, \texttt{label}, \texttt{acc}, and \texttt{nadd}. 
To run the code, type:
\begin{lstlisting}[style=notebook]
    ./parallel_zeta.sh operator_name primes scaling label acc nadd
\end{lstlisting}
Note that the last two options need not be included; they will take on default values in this case. However, if one wishes to include \texttt{nadd}, one must also include \texttt{acc}. If you would like to time your run, include \texttt{time} at the beginning of the above terminal input. The six options are explained below.

\begin{optDefnN}[label=newNumber]{\texttt{operator\_name}}
	The identifier of the Picard--Fuchs operator in question, stored in the input file \texttt{input\_data.txt}. Each line of \texttt{input\_data.txt} is stored as:
    \begin{center}
        \texttt{operator\_name other\_data},
    \end{center}
    and by referring to \texttt{operator\_name}, the shell script knows which data to read in. For example, in the list of provided example operators, \texttt{4.1.1} corresponds to the mirror quintic\footnote{This is the AESZ \texttt{newnumber}. The AESZ database can be found on the CYCluster website \cite{CYcluster}.}. A \texttt{jupyter} notebook (\texttt{example.ipynb}) has been provided in the \texttt{scripts} folder so that one may append one's own operators to the file \texttt{input\_data.txt} in the requisite format.
\end{optDefnN}

\begin{optDefnN}[label=primes]{\texttt{primes}}
	The range of primes for which to compute the Euler factors $R^{(b-1)}_p(\cL,T)$ for each $\varphi\in\IF_p^*$. This is input as a range\footnote{Be aware that in some terminals the square brackets must be input as \texttt{$\backslash$}\texttt{[} \texttt{$\backslash$}\texttt{]}.} 
    \texttt{[n\_min, n\_max]} and the script will compute $R^{(b-1)}_p(\cL,T)$ for all primes in the range $p_{n_{\text{min}}} \leq p \leq p_{n_{\text{max}}}$, where $p_n$ denotes the $n$'th prime.
\end{optDefnN}

\begin{optDefnN}[label=scaling]{\texttt{scaling}}
	The rational number $C$ such that the series $f_i(\co)$ are computed to $\lceil C p\rceil + \text{\texttt{nadd}}$ terms (see eq.~\eqref{eq:fWithC}). Should be input as a decimal, e.g. 1.5.
\end{optDefnN}

\begin{optDefnN}[label=id]{\texttt{label}}
	A user-defined label for the operator in question. 
    This determines the names of the files that are output by the bash script.
    For example, the coefficients of $R^{(b-1)}_p(\cL,T)$ for the options specified above will be output to the file \texttt{outputs/outputs\_label.txt}.
\end{optDefnN}

\begin{optDefnN}[label=intermedAcc]{\texttt{acc}}
    (Optional) The initial accuracy $A$ for the $p$-adically truncated recurrence (see \Cref{sec:pAdicTrunc}). If this is not included, the default value from eq.~\eqref{eq:A_lower_bound} is used in the computation.
\end{optDefnN}

\begin{optDefnN}[label=nadd]{\texttt{nadd}}
    (Optional) The number of terms that should vanish at the end of a series for the purposes of checking that the series terminates. If this is nonzero, then the default accuracy $A$ computed in eq.~\eqref{eq:A_lower_bound} may not be sufficient to guarantee accurate result, and so \texttt{acc} must be set appropriately.
\end{optDefnN}

After running the bash script from the \texttt{scripts} folder, the output file will be found in the \texttt{outputs} subfolder, saved as $\texttt{outputs\_label.txt}$, where \texttt{label} was Option \ref{id} input into the command. The outputs file contains the output of the main function \texttt{L\_functions}, where line-by-line, one finds the tuples $[p,\co_*,[a_p^{(1)}(\co_*),a_p^{(2)}(\co_*)]]$ for different primes $p$ and values of $\co_* \in \mathbb{F}_p^*$.

As future support for larger classes of operators of Calabi--Yau type is planned, an up-to-date version of the documentation will be kept on the \texttt{github} page. We also include an example \texttt{jupyter} notebook, \texttt{example.ipynb}, to demonstrate how to properly set up the code and convert a Calabi--Yau-type operator into the requisite form. This is kept in the \texttt{scripts} folder.

To conclude, we note that this package assumes a sufficiently nice form of the $\mtU$-matrix denominator, as in \Cref{foot:U_denominator}. Furthermore, as throughout \Cref{sec:pAdicTrunc}, we assume $p>C$. Thirdly, the code currently assumes $\ord(W(\co)^{-1})=0$; this is true for all but finitely many primes, and is a work-in-progress to include more generally. Finally, the \texttt{scripts} folder contains a \texttt{logs} subfolder, in which every time the main function \texttt{L\_functions} is run, a tuple $[p,\texttt{trunc\_deg},M]$ containing the prime $p$, the maximal of the degree of the entries of the matrix $\cU_p(\co)$ defined in eq.~\ref{eq:UratMatrix}, and the degree $M$ in $\co$ to which the series $f_i^{[M]}(\co)$ were computed (i.e. $M = \lceil Cp \rceil + \texttt{nadd}$). Note that if this degree is not $\leq \lceil C p \rceil$, a warning is printed.

\newpage

\bibliographystyle{alphaurl}
\bibliography{ref}
\vfill

\end{document}